\documentclass[11pt]{amsart}

\usepackage{amsfonts}
\usepackage{amsmath}
\usepackage{amssymb}
\usepackage{amsthm}
\usepackage{amsxtra}
\usepackage{epic, eepic}
\usepackage{epsfig,color}
\usepackage{fullpage}
\usepackage{vmargin}
\theoremstyle{plain}
\newtheorem{Thm}{Theorem}[section]

\newtheorem{Lem}[Thm]{Lemma}

\newtheorem{Conj}{Conjecture}

\theoremstyle{definition}

\theoremstyle{remark}
\newtheorem{Rem}{Remark}
\newtheorem*{Remrk}{Remark}

\setmarginsrb{33mm}{30mm}{33mm}{40mm}%
            {0mm}{10mm}{0mm}{10mm}

\begin{document}

\title{On a conjecture of Hivert and Thi\'ery about Steenrod operators}

\author{Michele D'Adderio}\author{Luca Moci}
\address{Max-Planck-Institut f\"{u}r Mathematik\\
Vivatsgasse 7, 53111 Bonn\\
Germany}\email{mdadderio@yahoo.it}

\address{Institut f\"{u}r Mathematik\\ Technische Universit\"{a}t
Berlin\\
Stra$\ss$e des 17. Juni, 136, Berlin\\
Germany}\email{moci@math.tu-berlin.de}

\begin{abstract}
We prove some results related to a conjecture of Hivert and
Thi\'ery about the dimension of the space of $q$-harmonics
(\cite{hivert}). In the process we compute the actions of the
involved operators on symmetric and alternating functions, which
have some independent interest. We then use these computations to
prove other results related to the same conjecture.
\end{abstract}

\maketitle


\section{Introduction}

The so called \textit{harmonics polynomials} (or
$\frak{S}_n$-\textit{harmonics}) are a classical object in
invariant and representation theory. They are the polynomial
solutions to the system of partial differential equations
$$
\nabla_k f(\mathbf{x}) =0\quad \text{for }k\geq 1,
$$
where $\mathbf{x}=x_1,x_2,\dots,x_n$ and the operators
$$
\nabla_k:=\sum_{i=1}^{n}\frac{\partial^k}{\partial x_{i}^k}
$$
are generalized laplacians. Since the $\nabla_k$'s are symmetric,
we have an action of the symmetric group $\frak{S}_n$ by
permutation of the variables. Hence the space of harmonic
polynomials is a representation of $\frak{S}_n$, that turns out to
be a regular representation, whose Frobenius characteristic is
(see \cite{macdonald})
$$
F_{n}(t)=F_{n;0}(t)=\sum_{\lambda \vdash n}s_{\lambda}\sum_{T\in
ST(n) }t^{co(T)},
$$
where $\lambda \vdash n$ indicates that $\lambda$ is a partition
of $n$, $s_{\lambda}$ is the Schur function indexed by $\lambda$,
$ST(\lambda)$ denotes the set of standard tableaux of shape
$\lambda$, and $co(T)$ denote the cocharge of the tableau $T$.

Recently many authors have studied various generalizations of the
operators $\nabla_k$'s, looking at similar spaces of polynomials.
It turns out that in many situations these spaces have
conjecturally the same Hilbert series (or the Frobenius
characteristic when the operators are symmetric) of the classical
harmonic polynomials.

In \cite{wood1,wood2,wood3} Wood raised several questions about
the \textit{rational Steenrod algebra} (twisted by the algebraic
Thom map), which is the subalgebra of the Weyl algebra generated
by the \textit{Steenrod operators}
$$
D^{*}_{k}=\sum_{i=1}^{n}x_{i}^{k}\left(1+x_i\frac{\partial}{\partial
x_i}\right),
$$
for $k\geq 1$. Let's call \textit{$1$-harmonic polynomials} the
ones killed by the duals of the $D^{*}_{k}$'s with respect to the
scalar product defined by
$$
\langle f(\mathbf{x}),g(\mathbf{x})\rangle:=f(\partial)
g(\mathbf{x})\big{|}_{\mathbf{x}=0},
$$
where $f(\partial)$ denote the differential operator obtained from
$f(\mathbf{x})$ by substituting the variables $x_i$ with the
operators $\frac{\partial}{\partial x_i}$. Among other things,
Wood asked (in a different language) if the space of $1$-harmonic
polynomials is a graded regular representation of the symmetric
group $\frak{S}_n$ (\textit{Rational hit conjecture}). We refer to
the works of Wood for motivations in Algebraic Topology.

In \cite{hivert} Hivert and Thi\'ery considered a deformed version
of those operators (and their duals), introducing the
\textit{$q$-Steenrod algebra}. They investigated questions similar
to the ones that Wood asked, finding interesting phenomena:
consider the operators
$$
D_{k;q}:=q\widetilde{D}_k+\nabla_k,
$$
with $\widetilde{D}_k:=\sum_{i=1}^n\, x_i
\partial_{i}^{k+1}$ and
$\nabla_k:=\sum_{i=1}^n\partial_{i}^{k}$, where
$\partial_j:=\frac{\partial}{\partial x_j}$, acting on
$\mathbb{C}(q)[\mathbf{x}]:=\mathbb{C}(q)[x_1,\dots,x_n]$, and $q$
is an indeterminate or a complex number.

We put
$$
\mathcal{H}_{\mathbf{x};q}:=\{g\in \mathbb{C}(q)[\mathbf{x}] \mid
D_{k;q}f=0\text{ for all $k\geq 1$}\},
$$
and we call its elements \textit{$q$-harmonics}. Also, we denote by
$$
\sum_{d\geq 0}\dim  \pi_d(\mathcal{H}_{\mathbf{x};q})t^d
$$
its Hilbert series.

Notice that the group $\frak{S}_n$ acts on these spaces by
permutation of the variables, since the operators involved are
symmetric.

\begin{Remrk}
Observe that for $q=0$ we retrieve the $\frak{S}_n$-harmonics,
while for $q=1$ the $D_{k;1}$'s are the dual of the Steenrod
operators. In fact the idea of Hivert and Thi\'ery was to
``interpolate'' the two situations via the coefficient $q$.
\end{Remrk}

In \cite{hivert} Hivert and Thi\'ery proved the following theorem
and stated the following conjecture.
\begin{Thm}[\cite{hivert}] \label{thmhivert}
When $q$ in an indeterminate, if we denote by $[n]_t !$ the usual
$t$-analogue of $n$-factorial, we have
$$
\sum_{d\geq 0}\dim  \pi_d(\mathcal{H}_{\mathbf{x};q})t^d << [n]_t
!
$$
with '$<<$' denoting coefficient-wise inequality.
\end{Thm}
In fact from this theorem it follows (see \cite{garsia}) that in
this case $\mathcal{H}_{\mathbf{x};q}$ is isomorphic to a graded
$\frak{S}_n$-submodule of the $\frak{S}_n$-harmonics.

\begin{Conj} \label{conj}
In the case where $q$ is a variable or a complex number not of the
form $-a/b$ where $a\in \{1,2,\dots ,n\}$ and $b\in \mathbb{N}$,
we have the equality
$$
\sum_{d\geq 0}\dim  \pi_d(\mathcal{H}_{\mathbf{x};q})t^d = [n]_t
!.
$$
In particular, in the case where $q$ is a variable,
$\mathcal{H}_{\mathbf{x};q}$ is isomorphic as a graded
$\frak{S}_n$-module to the $\frak{S}_n$-harmonics.
\end{Conj}
Notice that in the case where $q$ is a complex number the same
inequality of Theorem \ref{thmhivert} is not known even for
generic values of $q$.

After this work, in \cite{garsia} Bergeron, Garsia and Wallach
investigated even more general operators, bringing new insights in
this subject. Among other things, using commutative algebra, they
proved the following theorem.
\begin{Thm}[\cite{garsia}]
For any value of $q\in \mathbb{C}$ the dimension of the space of
$q$-harmonics in $n$ variables does not exceed $(n+1)!$.
\end{Thm}
Notice that of course the conjectured dimension for generic values
of $q\in \mathbb{C}$ is $n!$.

The common feature of all these works is the appearance of a
graded representations of $\frak{S}_n$ which is conjecturally
isomorphic to the classical $\frak{S}_n$-harmonics.

The present work arose from an attempt to make some progress on
Conjecture 1.

\subsection{First reductions}

\textit{Unless otherwise stated, $q$ will always be an
indeterminate. We will discuss the case $q\in \mathbb{C}$ mainly
in the last section of the present work.}\newline

We start with a general remark. Let $f\in
\mathcal{H}_{\mathbf{x};q}$. By multiplying by an element of
$\mathbb{C}(q)$, we can always assume that
$$
f=f_0+f_1q+f_2q^2+\cdots +f_mq^m
$$
with $f_i\in \mathbb{C}[\mathbf{x}]$ for all $i=1,\dots,m$ and
$f_0\neq 0\neq f_m$. It's easy to see (cf. \cite{garsia} or see
later) that $f_0$ is necessarily an $\frak{S}_n$-harmonic.

\begin{Lem}
Conjecture 1 is true if and only if for any $\frak{S}_n$-harmonic
$g$ we have a $q$-harmonic $f$ with $f_0=g$.
\end{Lem}
\begin{proof}
Suppose that the conjecture is true, and fix a basis
$g_{1},\dots,g_{n!}$. We can assume that each $g_i$ is of the form
$$
g_i=\sum_{j=0}^{m_i}g_{i,j}q^j,
$$
with $g_{i,j}\in \mathbb{C}[\mathbf{x}]$ and $g_{i,0}\neq 0 \neq
g_{i,m_i}$ for all $i$. We can also assume that the sequence
$\underline{m}:=(m_1,m_2,\dots,m_{n!})$ is in increasing order.
Choose a basis with minimal $\underline{m}$ with respect to the
lexicographic order. We claim that
$\{g_{1,0},g_{2,0},\dots,g_{n!,0}\}$ is a basis for the
$\frak{S}_n$-harmonics. If not, then we can find a non-trivial
linear combination
$$
\sum_{i=1}^{n!}\alpha_i g_{i,0}=0.
$$
But then we can replace $g_{n!}$ by the linear combination
$$
\sum_{i=1}^{n!}\alpha_i g_{i},
$$
and after dividing by a suitable power of $q$ we get a new basis,
with a smaller $\underline{m}$, which gives a contradiction. From
this the ``only if'' part follows.

The other implication is similar: choose a basis
$\{g_{1,0},g_{2,0},\dots,g_{n!,0}\}$ of the
$\frak{S}_n$-harmonics, and by using the hypothesis we can find
$q$-harmonics $g_{1},\dots,g_{n!}$ such that
$$
g_i=\sum_{j=0}^{m_i}g_{i,j}q^j.
$$
I claim that these are independent over $\mathbb{C}(q)$. If not,
we would have a nontrivial combination
$$
\sum_{i=1}^{n!}\alpha_i(q) g_{i}=0,
$$
with $\alpha_i(q) \in \mathbb{C}(q)$. Of course we can normalize
these coefficients so that they are all polynomials, and at least
one non-zero coefficient has non-zero constant term. But then the
constant term of this linear combination would give a non-trivial
linear relation among the $g_{i,0}$'s, which gives a
contradiction.
\end{proof}

From the easy relations
$$
[D_{k;q},D_{h;q}]=q(k-h)D_{k+h;q},
$$
it follows that a polynomial $f$ is in
$\mathcal{H}_{\mathbf{x};q}$ if and only if
$$
D_{1;q}f=D_{2;q}f=0.
$$
This is clearly true even for $q\in \mathbb{C}$, $q\neq 0$.

It's easy to show (see \cite{garsia}) that the previous two
equations are equivalent to the following system of equations:
\begin{eqnarray}
\notag \nabla_k f_0 & = & 0,\\
\nabla_k f_i & = & -\widetilde{D}_k f_{i-1}\quad \text{for $i=1,2,\dots,m$},\\
\notag \widetilde{D}_k f_m & = & 0,
\end{eqnarray}
for $k=1,2$. Notice in particular that $f_0$ is an
$\frak{S}_n$-harmonic.

Together with the previous lemma, this shows that if for any
$\frak{S}_n$-harmonic $f_0$ we are able to find $f_1,f_2,\dots$
that satisfy those equations, then the conjecture is true.

In this work we try to attack Conjecture 1 using these
observations. The idea would be to construct the entire sequence
$f_1,f_2,\dots$ for any $f_0$. We only succeeded in constructing
an $f_1$ for any $f_0$, and the corresponding $f_2$ for some
special $\frak{S}_n$-harmonic. We found two methods to achieve
this, one computationally heavier than the other, that provide
different solutions. We present both of them, since the hope is to
eventually find the entire sequence $f_1,f_2,\dots$.

Along the way we determine the action of the operators $\nabla_1$,
$\nabla_2$, $\widetilde{D}_1$ and $\widetilde{D}_2$ on symmetric
and alternating polynomials, which is of independent interest.

In fact in the last part we will use these actions to prove some
results related to the conjecture in the case $q\in \mathbb{C}$.

\subsection{Further reductions}

The first goal is to show how to construct an $f_1$ for any
$\frak{S}_n$-harmonic. Before doing that we want to show that it's
enough to construct an $f_1$ for $f_0=\partial_1\Delta$, where
$\Delta$ denotes the Vandermonde determinant in the variables
$x_1,\dots,x_n$.
\begin{Rem} \label{rmkvander}
In what follows we will repeatedly use the observation that any
symmetric homogeneous differential operators that lower the degree
kills the Vandermonde determinant. This is true since when we act
on $\Delta$ with such an operator we still get an alternant, but
of a lower degree. This forces it to be zero since the Vandermonde
determinant is the alternant of smallest possible degree.
\end{Rem}

Suppose that we know how to construct such an $f_1$. By permuting
its variables, it's clear how to construct an $f_1$ for
$\partial_i\Delta$ for all $i$'s. Let's call it $f_{1}^{(i)}$.

Also, remember that the partial derivatives of $\Delta$ span the
space of $\frak{S}_n$-harmonics. Hence by linearity it's enough to find an $f_1$ for any of those derivatives.

We set for any multi-index $\alpha=(\alpha_1,\alpha_2,\dots,\alpha_n)\in \mathbb{N}^n$
$$
\partial^{\alpha}:=\partial_{1}^{\alpha_1}\partial_{2}^{\alpha_2}\dots \partial_{n}^{\alpha_n}.
$$
\begin{Rem} \label{reductionrmk}
We have
$$
[\widetilde{D}_k,\partial^{\alpha}]=-\sum_{\alpha_i\neq 0}\alpha_i\partial^{\alpha+kv_i},
$$
where $v_i\in \mathbb{N}^n$ is the vector with $1$ in the $i$-th
position and $0$ elsewhere. Since clearly $\tilde{D}_k \Delta =0$
(cf. Remark \ref{rmkvander}), it follows that
$$
-\widetilde{D}_k\partial^{\alpha}\Delta=\sum_{\alpha_i\neq 0}\alpha_i\partial^{\alpha+kv_i}\Delta=\sum_{\alpha_i\neq 0}\alpha_i\partial^{\alpha-v_i}(\partial_{i}^{k+1}\Delta).
$$
Hence if we set
$$
f_{1}^{\alpha}:=\sum_{\alpha_i\neq 0}\alpha_i\partial^{\alpha-v_i}f_{1}^{(i)},
$$
we have
$$
\nabla_k f_{1}^{\alpha}=-\widetilde{D}_k\partial^{\alpha}\Delta
$$
for all multi-indices $\alpha$ and $k=1,2$.
\end{Rem}

We are left with the task of computing $f_{1}^{(1)}$.

\subsection{Organization of the paper}

The rest of the paper is organized in the following way:
\begin{itemize}
    \item In the second section we find an $f_1$ for $\partial_1 \Delta$.
    \item In the third section we find an entire family of $f_1$'s, which
    include the previous one as a special case. For one member of this family we
    find an $f_2$ also, but we relegated the computations in the appendix.
    \item In the fourth section we show another method of finding an
    $f_1$ and an $f_2$ for $\partial_1 \Delta$.
    \item In the fifth section we compute systematically the
    action of the operators $\nabla_1$, $\nabla_2$,
    $\widetilde{D}_1$ and $\widetilde{D}_2$ on symmetric and
    alternating polynomials.
    \item In the sixth section we discuss the case $q\in \mathbb{C}$.
    We apply our formulae to investigate what we will call ``singular''
    values of $q$. We prove that most of the values excluded in
    Conjecture 1 are indeed singular, and we finally state a new
    conjecture on these singular values.
\end{itemize}

\section{Computation of $f_1$ for $\partial_1\Delta$}

We want to construct an $f_1=f_{1}^{(1)}$ for $f_0=\partial_1\Delta$. We can of course assume that $f_1$ is homogeneous.

We want
$$
\nabla_k f_1=-\tilde{D}_k \partial_1 \Delta
$$
for $k=1,2$. We already noticed that
$[\tilde{D}_k,\partial_1]=-\partial_{1}^{k+1}$, so we can rewrite those equations as
$$
\nabla_k f_1=\partial_{1}^{k+1} \Delta.
$$
We now assume that $f_1$ is of the form $\Delta^{(1)}g$, where
$\Delta^{(1)}$ is the Vandermonde in the variables $x_2,\dots,x_n$
and $g$ is a polynomial of the form
$$
g=\sum_{j=1}^{n-2}g_j x_{1}^j
$$
where each $g_j$ is a symmetric polynomial in $x_2,\dots,x_n$
homogeneous of degree $n-2-j$. In this case we get
\begin{eqnarray*}
\nabla_1 f_1 & = & (\nabla_1 \Delta^{(1)})g+\Delta^{(1)}(\nabla_1 g) \\
 & = & \Delta^{(1)}\left(\sum_{s=0}^{n-2}(\nabla_1 g_s )x_{1}^s + \sum_{j=1}^{n-2}j g_j x_{1}^{j-1} \right)\\
 & = & \Delta^{(1)}\left(\sum_{s=0}^{n-3}\left(\nabla_1 g_s  + (s+1) g_{s+1}
 \right)x_{1}^{s}\right)
\end{eqnarray*}
where the second equality holds since $\nabla_1
\Delta^{(1)}=0$.\newline

We fix the notation $e_k:=e_k(x_2,\dots,x_n)$, which will be used
also in the following sections, except the last one. We start by
recording some easy identities:
\begin{eqnarray} \label{nabla}
\notag \nabla_1 e_k & = & (n-k)e_{k-1};\\
\nabla_{1}^s e_k & = & (n-k)(n-k+1)\cdots (n-k+s-1)e_{k-s}\quad \text{for $s\geq 1$};\\
\notag \nabla_1 e_k^a & = & a(n-k)e_{k}^{a-1}e_{k-1};\\
\notag \nabla_h e_k& = & 0\qquad\text{for all $h\geq 2$.}
\end{eqnarray}

We have
\begin{eqnarray*}
\partial_{1}^2 \Delta & = & \Delta^{(1)}\left(\partial_{1}^2\prod_{j=2}^{n}(x_1-x_j)\right) \\
 & = & \Delta^{(1)}\left(\partial_{1}^2\sum_{j=0}^{n-1}(-1)^{j}e_jx_{1}^{n-1-j}\right) \\
 & = & \Delta^{(1)}\left(\sum_{j=0}^{n-3}(-1)^{j}(n-1-j)(n-2-j)e_jx_{1}^{n-3-j}\right)\\
 & = & \Delta^{(1)}\left(\sum_{s=0}^{n-3}(-1)^{n-3-s}(s+2)(s+1)e_{n-3-s}x_{1}^{s}\right).
\end{eqnarray*}
Equating the coefficients we get the system of equations
\begin{Lem}
\begin{equation} \label{C1}
\tag{\textbf{C}1} (-1)^{n-3-s}(s+2)(s+1)e_{n-3-s}=\nabla_1 g_s +
(s+1)g_{s+1}\quad\text{for $s=0,1,\dots,n-3$.}
\end{equation}
\end{Lem}
This system can be integrated in many ways. We now use these
equations to write all the $g_j$'s for $j\geq 1$ in terms of
$\nabla_{1}^h g_0$ for $h\geq 0$.
\begin{Lem}
For $s=1,2,\dots,n-2$ we have the following formula:
\begin{equation} \label{gss}
\tag{$\bullet$} g_s=(-1)^{n-2-s}(s+1)s\,
e_{n-2-s}+\frac{(-1)^{s}}{s!}\nabla_{1}^{s}g_0.
\end{equation}
\end{Lem}
\begin{proof}
First of all notice that for $s\geq 0$ we can write the equations
(\ref{C1}) as
$$
g_{s+1}=(-1)^{n-3-s}(s+2) e_{n-3-s}-\frac{1}{s+1}\nabla_1 g_s.
$$
We proceed by induction on $s$, the case $s=1$ being just equation
(\ref{C1}). Assume that the result is true for $s\geq 1$. Then we
have
\begin{eqnarray*}
g_{s+1} & = & (-1)^{n-3-s}(s+2) e_{n-3-s}-\frac{1}{s+1}\nabla_1 g_s \\
 & = &(-1)^{n-3-s}(s+2)
 e_{n-3-s}-\frac{1}{s+1}\nabla_1\left((-1)^{n-2-s}(s+1)s\,
e_{n-2-s}+\frac{(-1)^s}{s!}\nabla_{1}^{s}g_0\right) \\
 & = & (-1)^{n-3-s}((s+2)+s(s+2))\,
e_{n-3-s}+\frac{(-1)^{s+1}}{(s+1)!}\nabla_{1}^{s+1}g_0\\
 & = & (-1)^{n-3-s}(s+2)(s+1)\,
e_{n-3-s}+\frac{(-1)^{s+1}}{(s+1)!}\nabla_{1}^{s+1}g_0.
\end{eqnarray*}

\end{proof}

Of course to find what we want, we need to take into account the
other set of equations coming from $\nabla_2
f_1=\partial_{1}^3\Delta$. We have
\begin{eqnarray*}
\sum_{i=1}^{n}\partial_{i}^2 f_1 & = & \sum_{i=1}^{n}\left(
(\partial_{i}^2 \Delta^{(1)})g+\Delta^{(1)}
(\partial_{i}^2 g)+ 2(\partial_{i} \Delta^{(1)}\cdot \partial_{i} g)\right) \\
 & = & (\nabla_2 \Delta^{(1)})g+\Delta^{(1)}
(\nabla_2 g)+ \sum_{i=1}^{n}2(\partial_{i} \Delta^{(1)}\cdot \partial_{i} g)\\
 & = & \Delta^{(1)}(\nabla_2 g)+ \sum_{i=2}^{n}2(\partial_{i} \Delta^{(1)}\cdot
\partial_{i} g).
\end{eqnarray*}
Dividing by $\Delta^{(1)}$ we get
\begin{eqnarray*}
\frac{1}{\Delta^{(1)}}\sum_{i=1}^{n}\partial_{i}^2 f_1 & = &
\nabla_2 g+ 2\sum_{i=2}^{n}\left(\frac{\partial_{i}
\Delta^{(1)}}{\Delta^{(1)}}\cdot
\partial_{i} g\right)\\
 & = &
\nabla_2 g+ 2\sum_{i=2}^{n}\left(\partial_{i}\log
\Delta^{(1)}\cdot
\partial_{i} g\right)\\
 & = &
\nabla_2 g+ 2\sum_{i=2}^{n}\left(\,\,\sum_{2\leq j\leq n,\,\,
j\neq i}\frac{(-1)^{\chi (i>j)}}{x_i-x_j}\cdot
\partial_{i} g\right)\\
 & = & \nabla_2 g+ 2\sum_{2\leq i<j\leq n}\frac{1}{x_i-x_j}\left(
 \partial_{i}-\partial_{j}\right)g,
\end{eqnarray*}
where $\chi(\mathcal{P})$ is equal to $1$ if the proposition
$\mathcal{P}$ is true, $0$ otherwise.

Setting
$$
P_2:=\sum_{2\leq i<j\leq n}\frac{1}{x_i-x_j}\left(
\partial_{i}-\partial_{j}\right),
$$
we have
\begin{eqnarray*}
\frac{1}{\Delta^{(1)}}\sum_{i=1}^{n}\partial_{i}^2 f_1 & = &
(\nabla_2 + 2P_2) g\\
 & = & \sum_{s=0}^{n-4}\left((\nabla_2+2P_2)g_s+(s+2)(s+1)g_{s+2}\right)x_{1}^s.
\end{eqnarray*}

On the other hand we have
\begin{eqnarray*}
\partial_{1}^3 \Delta & = & \Delta^{(1)}\left(\partial_{1}^3\prod_{j=2}^{n}(x_1-x_j)\right) \\
 & = & \Delta^{(1)}\left(\partial_{1}^3\sum_{j=0}^{n-1}(-1)^{j}e_jx_{1}^{n-1-j}\right) \\
 & = & \Delta^{(1)}\cdot \sum_{j=0}^{n-4}(-1)^{j}(n-1-j)(n-2-j)(n-3-j)e_jx_{1}^{n-4-j}\\
 & = & \Delta^{(1)}\cdot\sum_{s=0}^{n-4}(-1)^{n-4-s}(s+3)(s+2)(s+1)e_{n-4-s}\,x_{1}^{s}.
\end{eqnarray*}
Equating the coefficients we get the following system of
equalities:
\begin{Lem}
\begin{equation} \label{C2}
\tag{\textbf{C}2}(-1)^{n-4-s}(s+3)(s+2)(s+1)e_{n-4-s} =
(\nabla_2+2P_2)g_s+(s+2)(s+1)g_{s+2},
\end{equation}
$$
\,\,\,\qquad\text{for
$s=0,1,\dots,n-4$.}\qquad\qquad\qquad\qquad\qquad\qquad\qquad\qquad\qquad\qquad\qquad
$$
\end{Lem}
We study now some properties of the operator $P_2$.
\begin{Lem}
We have the following identities:
\begin{eqnarray} \label{p2}
P_2 e_k & = & -\left(\!\!\!\begin{array}{c}
 n-k+1 \\
  2 \\
\end{array} \!\!\!\right) e_{k-2}; \\
\notag P_2 e_{k}^h & = & he_{k}^{h-1}P_2
e_k=-h\left(\!\!\!\begin{array}{c}
 n-k+1 \\
  2 \\
\end{array} \!\!\!\right) e_{k}^{h-1} e_{k-2}.
\end{eqnarray}
\end{Lem}
\begin{proof}
If we denote by $e_{k}^{(i)}$ the elementary symmetric function of
degree $k$ in the variables $\{x_2,\dots,x_n\}\setminus \{x_i\}$,
we have
$$
\partial_ie_k=e_{k-1}^{(i)}.
$$
Consider the difference
$$
\partial_i e_k-\partial_j e_k=
e_{k-1}^{(i)}-e_{k-1}^{(j)}.
$$
The monomials in $e_{k-1}^{(i)}$ that don't involve $x_j$ are
cancelled by the ones in $e_{k-1}^{(j)}$ that don't contain $i$;
while the monomials in $e_{k-1}^{(i)}$ that involve $x_j$ can be
paired with the ones in $e_{k-1}^{(j)}$ that involve $x_i$, to get
a factor $x_j-x_i$, so that when we divide by $x_i-x_j$ we are
left only with the negative of a multiple of $e_{k-2}$.

To see what this multiple is, it's enough to count how many times
the monomial $x_2x_3\cdots x_{k-1}$ appears: this number is the
number of ways of choosing $i$ and $j$ in $\{k,k+1,\dots,n\}$,
which is what we wanted.

The second identity follows from the first one and Leibniz rule.
\end{proof}
\begin{Lem}
If $g$ is a symmetric polynomial, then
$$
[\nabla_1,P_2]g=0.
$$
\end{Lem}
\begin{proof}
It's enough to check this relation on the monomials $e_{\lambda}$,
where $\lambda$ denotes as usual a partition, since they form a
basis of symmetric polynomials. Using repeatedly Leibniz rule we
reduce ourselves to check the identity on the $e_k$'s. But this
follows immediately from the identities (\ref{nabla}) and
(\ref{p2}).
\end{proof}

Substituting (\ref{gss}) in (\ref{C2}) and using the previous lemmas we get
$$
(-1)^{n-4-s}(s+3)(s+2)(s+1)e_{n-4-s}=
\qquad\qquad\qquad\qquad\qquad\qquad\qquad\qquad\qquad\qquad\qquad\qquad
$$
\begin{eqnarray*}
\qquad\qquad & = & (\nabla_2+2P_2)\left((-1)^{n-2-s}(s+1)s\,
e_{n-2-s}+\frac{(-1)^s}{s!}\nabla_{1}^{s}g_0\right)\\ & + &
(s+2)(s+1)\left((-1)^{n-4-s}(s+3)(s+2)\,
e_{n-4-s}+\frac{(-1)^{s+2}}{(s+2)!}\nabla_{1}^{s+2}g_0\right)\\
 & = & 2(-1)^{n-3-s}(s+1)s\left(\!\!\!\begin{array}{c}
 s+3 \\
  2 \\
\end{array} \!\!\!\right)e_{n-4+s}+ \frac{(-1)^s}{s!}(\nabla_2+2P_2)\nabla_{1}^{s}g_0\\
 & + & (-1)^{n-4-s}(s+3)(s+2)^2(s+1)\,
e_{n-4-s}+ \frac{(-1)^s}{s!}\nabla_{1}^{s+2}g_0
\end{eqnarray*}
\begin{eqnarray*}
 & = & (-1)^{n-3-s}(s+3)(s+2)(s+1)s\, e_{n-4+s}+ \frac{(-1)^s}{s!}(\nabla_2+2P_2)\nabla_{1}^{s}g_0\\
 & + & (-1)^{n-4-s}(s+3)(s+2)^2(s+1)\,
e_{n-4-s}+ \frac{(-1)^s}{s!}\nabla_{1}^{s+2}g_0\\
 & = & 2(-1)^{n-4-s}(s+3)(s+2)(s+1) e_{n-4-s}+
 \frac{(-1)^s}{s!}\nabla_{1}^{s}(\nabla_2+2P_2+\nabla_{1}^{2})g_0,
\end{eqnarray*}
from which we get the following system of identities:
$$
(-1)^{n-4-s}(s+3)(s+2)(s+1)\, e_{n-4-s}+
 \frac{(-1)^s}{s!}\nabla_{1}^{s}(\nabla_2+2P_2+\nabla_{1}^{2})g_0=0
$$
$$
\text{for
$s=0,1,\dots,n-4$.}\qquad\qquad\qquad\qquad\qquad\qquad\qquad\qquad\qquad\qquad\qquad
$$

These equations can be rewritten in the following form:
\begin{Lem}
$$
\nabla_{1}^{s}(\nabla_2+2P_2+\nabla_{1}^{2})g_0=(-1)^{n-1}(s+3)!\,
e_{n-4-s}
$$
$$
\text{for
$s=0,1,\dots,n-4$.}\qquad\qquad\qquad\qquad\qquad\qquad\qquad\qquad\qquad\qquad\qquad
$$
\end{Lem}
Notice that by (\ref{nabla}) we have
$$
\nabla_{1}^s e_{n-4}=4\cdot 5\cdot \cdots \cdot (s+3)
e_{n-4-s}=\frac{1}{6}(s+3)!e_{n-4-s},
$$
hence
$$
(\nabla_2+2P_2+\nabla_{1}^{2})g_0=(-1)^{n-1}6\, e_{n-4}
$$
would give a solution to all our systems.

\begin{Rem} \label{remarkonnabla1andP2}
It's straightforward to check that
$$
(\nabla_{1}^2+2P_2)e_{k}=0\quad \text{for all $k$.}
$$
\end{Rem}
Since also $\nabla_2 e_k=0$ for all $k$, we must look for a $g_0$
that involves $e_{\lambda}$ with partitions $\lambda$ consisting
of at least two parts.

In the following calculations we will use identities (\ref{nabla})
and (\ref{p2}); remember that the $e_k$'s are in the $n-1$
variables $x_2,\dots,x_n$.
$$
2P_2 (e_{n-3}e_1) = 2P_2 (e_{n-3})e_1+2e_{n-3}P_2(e_1) =-12
e_{n-5}e_1;
$$
\begin{eqnarray*}
\nabla_{1}^{2}(e_{n-3}e_1) & = & (\nabla_{1}^2 e_{n-3})e_1+2
\nabla_{1}e_{n-3}\nabla_{1}e_1+ e_{n-3}(\nabla_{1}^2e_1)\\
 & = & 12 e_{n-5}e_1+ 6(n-1) e_{n-4};
\end{eqnarray*}
\begin{eqnarray*}
\nabla_{2}(e_{n-3}e_1) & = & \sum_{i=2}^{n} (\partial_{i}^2
e_{n-3})e_1 + 2 \sum_{i=2}^{n}
\partial_{i} e_{n-3}\partial_{i} e_1+ \sum_{i=2}^{n} e_{n-3}(\partial_{i}^2
e_{1})\\
 & = & 2\nabla_1 e_{n-3}=6 e_{n-4}.
\end{eqnarray*}

From these we get
$$
(\nabla_2+2P_2+\nabla_{1}^{2})e_{n-3}e_1= 6 n e_{n-4}.
$$
Hence our solution will be
$$
g_0:=\frac{(-1)^{n-1}}{n}e_{n-3}e_1.
$$

Now we want to make formula (\ref{gss}) more explicit.
\begin{Lem}
For $s\geq 1$ we have
$$
\nabla_{1}^s(e_{n-3}e_1)=\frac{(s+2)!}{2}e_{n-3-s}e_1+\frac{(s+1)!}{2}s(n-1)e_{n-2-s}.
$$
\end{Lem}
\begin{proof}
By induction on $s$, the case $s=1$ being clear. We assume the
formula true for $s\geq 1$. We have
\begin{eqnarray*}
\nabla_{1}^{s+1}(e_{n-3}e_1) & = & \nabla_1\left(
\frac{(s+2)!}{2}e_{n-3-s}e_1+\frac{(s+1)!}{2}s(n-1)e_{n-2-s}\right)\\
 & = & \frac{(s+2)!}{2}((\nabla_1e_{n-3-s})e_1+e_{n-3-s}(\nabla_1 e_1))+\frac{(s+1)!}{2}s(n-1)\nabla_1 e_{n-2-s}\\
 & = &
 \frac{(s+3)!}{2}e_{n-4-s}e_1+\frac{(s+2)!}{2}(n-1)e_{n-3-s}+\frac{(s+2)!}{2}s(n-1)e_{n-3-s}\\
 & = &
 \frac{(s+3)!}{2}e_{n-4-s}e_1+\frac{(s+2)!}{2}(s+1)(n-1)e_{n-3-s}.
\end{eqnarray*}
\end{proof}
Plugging these formulae into (\ref{gss}) we get for all $s\geq 1$
\begin{eqnarray*} g_s & = & (-1)^{n-2-s}(s+1)s\,
e_{n-2-s}+\frac{(-1)^s}{s!}\nabla_{1}^{s}g_0\\
 & = & (-1)^{n-2-s}(s+1)s\,
e_{n-2-s}+\frac{(-1)^s}{s!}\nabla_{1}^{s}\left(\frac{(-1)^{n-1}}{n}e_{n-3}e_1\right)\\
 & = & (-1)^{n-2-s}(s+1)s\,
e_{n-2-s}+\\
 & + & \frac{(-1)^{n+s-1}}{s!\,
n}\left(\frac{(s+2)!}{2}e_{n-3-s}e_1+\frac{(s+1)!}{2}s(n-1)e_{n-2-s}\right)\\
 & = & \left(\frac{n+1}{2n}\right)(-1)^{n-2-s}(s+1)s\,
e_{n-2-s}+\left(\frac{1}{2n}\right)(-1)^{n-1-s}(s+2)(s+1)\,
e_{n-3-s}e_1\\
 & = & \frac{(-1)^{n-2-s}}{n}\left((n+1)\left(\!\!\!\begin{array}{c}
 s+1 \\
  2 \\
\end{array} \!\!\!\right)\,
e_{n-2-s}-\left(\!\!\!\begin{array}{c}
 s+2 \\
  2 \\
\end{array} \!\!\!\right)
e_{n-3-s}e_1\right).
\end{eqnarray*}
We follow the convention that the binomial ``$n$ choose $k$'' is
$0$ when $n<k$, hence this formula works for $s\geq 0$.

Putting everything together, we get the formula
$$
f_1=f_{1}^{(1)}=\Delta^{(1)}\sum_{s=0}^{n-2}\frac{(-1)^{n-2-s}}{n}\left((n+1)\left(\!\!\!\begin{array}{c}
 s+1 \\
  2 \\
\end{array} \!\!\!\right)\,
e_{n-2-s}-\left(\!\!\!\begin{array}{c}
 s+2 \\
  2 \\
\end{array} \!\!\!\right)
e_{n-3-s}e_1\right)x_{1}^{s}.
$$

Encouraged by this promising first step, we tried to pursue our
methods to compute an $f_2$ for our $f_1$. Notice that this $f_2$
would work only for $\partial_j \Delta$, and not for a general
$\frak{S}_n$-harmonic, since the other reduction that we did for
$f_1$ doesn't work for $f_2$.

With some patience and stamina we went trough our computations, to
finally realize that we couldn't find an $f_2$ for all values of
$n$ in this way. But not all efforts were lost: some of those
computations are now part of the fifth section!

Looking back at the work in the present section, we realized that
something more general could be done.

\section{A family of $f_1$'s for $\partial_1 \Delta$}

When we constructed our explicit $f_1$ we had to solve the system
of equations
$$
(\nabla_2+2P_2+\nabla_{1}^{2})g_0=(-1)^{n-1}6\, e_{n-4}.
$$
Of course the solution that we had found was not unique. In fact
there are infinitely many solutions to this system. In this
section we construct a whole family of solutions. Of course we are
going to use much of what we did in the last section.

We need the following identities:
\begin{Lem}
For $k\geq h$ we have
\begin{eqnarray*}
\nabla_2(e_ke_2) & = & 2(n-k)e_{k-1}e_{1}-2ke_k;\\
(\nabla_{1}^2+2P_2)e_ke_2 & = & 2(n-k)(n-2)e_{k-1}e_1;\\
(\nabla_{1}^2+2P_2+\nabla_2)e_ke_2 & = &
2(n-k)(n-1)e_{k-1}e_1-2ke_k;\\
(\nabla_{1}^2+2P_2+\nabla_2)e_ke_{1}^2 & = & 4n(n-k)e_{k-1}e_1+
2n(n-1)e_k.
\end{eqnarray*}
\end{Lem}
\begin{proof}
The first identity is a special case of a more general formula
that can be found in the fifth section with its proof. The second
one follows easily from remark (\ref{remarkonnabla1andP2}). The
third one follows from the previous two. The last one is a special
case of previous identities.
\end{proof}
We can now look for a solution of our system. We assume that $g_0$
is of the form
$$
g_0=a\, e_{n-4}e_2+b\,e_{n-4}e_{1}^2+c\, e_{n-3}e_1,
$$
where $a=a(n)$, $b=b(n)$ and $c=c(n)$ are indeterminate
coefficients.

We have
\begin{eqnarray*}
(\nabla_{1}^2+2P_2+\nabla_2)g_0 & = &
a(8(n-1)e_{n-5}e_1-2(n-4)e_{n-4})\\
 & + & b(16n\, e_{n-5}e_1 + 2n(n-1)e_{n-4})\\
 & + & c\, 6n\, e_{n-4},
\end{eqnarray*}
from which we get the two equations
\begin{eqnarray*}
a\, 8(n-1)+b\, 16n & = & 0;\\
-a\, 2(n-4)+b\, 2n(n-1)+ c\, 6n & = & (-1)^{n-1}6.
\end{eqnarray*}
Solving for $a$ and $b$ we get
$$
a=-\frac{6((-1)^{n-1}-c\, n)}{n^2-7},\qquad
b=\frac{3(n-1)((-1)^{n-1}-c\, n)}{n(n^2-7)},
$$
where $c$ can be any number. Hence we get the family of solutions
$$
g_{0;c}=-\frac{6((-1)^{n-1}-c\, n)}{n^2-7}
e_{n-4}e_2+\frac{3(n-1)((-1)^{n-1}-c\,
n)}{n(n^2-7)}e_{n-4}e_{1}^2+c\, e_{n-3}e_1.
$$
Observe that in the previous section we got $g_{0;c}$ for
$c=(-1)^{n-1}/n$.

We record the following two identities, which are just
consequences of the identities that we already established and
Leibniz rule:
\begin{Lem}
\begin{eqnarray*}
\nabla_{1}^{s}e_{n-4}e_2 & = & \frac{(s+3)!}{3!}e_{n-s-4}e_2
+s(n-2)\frac{(s+2)!}{3!}e_{n-s-3}e_1\\
 & + & \frac{s(s-1)}{2}(n-1)(n-2)\frac{(s+1)!}{3!}e_{n-s-2};\\
\nabla_{1}^{s}e_{n-4}e_{1}^2 & = &
\frac{(s+3)!}{3!}e_{n-s-4}e_{1}^2 +
2s(n-1)\frac{(s+2)!}{3!}e_{n-s-3}e_1\\
 & + & s(s-1)(n-1)^2\frac{(s+1)!}{3!}e_{n-s-2}.
\end{eqnarray*}
\end{Lem}

Hence we have
\begin{eqnarray*}
\nabla_{1}^{s}g_{0;c} & = & -\frac{6((-1)^{n-1}-c\,
n)}{n^2-7}\left(
\frac{(s+3)!}{3!}e_{n-s-4}e_2+s(n-2)\frac{(s+2)!}{3!}e_{n-s-3}e_1\right.\\
 & + & \left.\frac{s(s-1)}{2}(n-1)(n-2)\frac{(s+1)!}{3!}e_{n-s-2}\right)\\
 & + & \frac{3(n-1)((-1)^{n-1}-c\, n)}{n(n^2-7)}\left( \frac{(s+3)!}{3!}e_{n-s-4}e_{1}^2 +
2s(n-1)\frac{(s+2)!}{3!}e_{n-s-3}e_1 \right.\\
 & + & \left. s(s-1)(n-1)^2\frac{(s+1)!}{3!}e_{n-s-2}\right)\\
 & + &
 c\left(\frac{(s+2)!}{2}e_{n-3-s}e_1+\frac{(s+1)!}{2}s(n-1)e_{n-s-2}\right)\\
 & = & -\frac{6((-1)^{n-1}-c\,
n)}{n^2-7} \frac{(s+3)!}{3!}e_{n-s-4}e_2+\frac{3(n-1)((-1)^{n-1}-c\, n)}{n(n^2-7)}\frac{(s+3)!}{3!}e_{n-s-4}e_{1}^2\\
 & + & \left(-\frac{((-1)^{n-1}-c\,
n)}{n^2-7}s(n-2)+\frac{(n-1)((-1)^{n-1}-c\, n)}{n(n^2-7)}s(n-1)+
\frac{c}{2}\right)(s+2)!e_{n-s-3}e_1\\
 & + & \left(-\frac{((-1)^{n-1}-c\,
n)}{n^2-7}\frac{s(s-1)}{2}(n-1)(n-2)\right.\\
 & + & \left.\frac{(n-1)((-1)^{n-1}-c\, n)}{n(n^2-7)}\frac{s(s-1)}{2}(n-1)^2+c\frac{s(n-1)}{2}
 \right)(s+1)!e_{n-s-2}\\
  & = & -\frac{((-1)^{n-1}-c\,
n)}{n^2-7} (s+3)!e_{n-s-4}e_2+\frac{(n-1)((-1)^{n-1}-c\, n)}{2n(n^2-7)}(s+3)!e_{n-s-4}e_{1}^2\\
 & + & \left(\frac{s((-1)^{n-1}-c\,
n)}{n(n^2-7)}+ \frac{c}{2}\right)(s+2)!e_{n-s-3}e_1\\
 & + & \left(\frac{((-1)^{n-1}-c\,
n)}{n^2-7}
\frac{s(s-1)}{2}\frac{n-1}{n}+c\frac{s(n-1)}{2}\right)(s+1)!e_{n-s-2}.
\end{eqnarray*}
Finally for $s\geq 1$ we have
\begin{eqnarray*}
g_{s;c} & = & (-1)^{n-2-s}(s+1)s\, e_{n-s-2}+\frac{(-1)^s}{s!}\nabla_{1}^sg_{0;c}\\
 & = & -(-1)^{s}\frac{((-1)^{n-1}-c\,
n)}{n^2-7} (s+3)(s+2)(s+1)e_{n-s-4}e_2\\
 & + & (-1)^s\frac{(n-1)((-1)^{n-1}-c\, n)}{2n(n^2-7)}(s+3)(s+2)(s+1)e_{n-s-4}e_{1}^2\\
 & + & (-1)^s\left(\frac{s((-1)^{n-1}-c\,
n)}{n(n^2-7)}+ \frac{c}{2}\right)(s+2)(s+1)e_{n-s-3}e_1\\
 & + & (-1)^s\left(\frac{((-1)^{n-1}-c\,
n)}{n^2-7}
\frac{s(s-1)}{2}\frac{n-1}{n}+c\frac{s(n-1)}{2}+s(-1)^n\right)(s+1)e_{n-s-2}.
\end{eqnarray*}
From this we could write a formula for $f_{1;c}$.

At this point we looked for a value of $c$ for which we could find
an $f_2$. In the end we found exactly one for each value of $n$:
$$
c=(-1)^{n-1}\frac{2(2n^3-2n-3)}{3n(n-1)(n^2+n+2)}.
$$
We relegated the derivation of the value of $c$ and the
computation of the corresponding $f_2$ in the appendix, since the
calculation is quite long. Reading the appendix should make clear
that these methods can't be pushed much further without a
tremendous stamina.

In the next section we show instead a different method to get
other $f_1$'s.

\section{Another computation of $f_1$}

We want to find an $f_1$ for $\partial_j \Delta$. In fact we will
prove something more. First of all we make the following simple
observation: from the obvious $\nabla_1 \Delta=0$ (see Remark
\ref{rmkvander}) we get
$$
\nabla_{1}^{(j)}\Delta=-\partial_j \Delta,
$$
where $\nabla_{1}^{(j)}$ denotes the sum of the partial
derivatives with $\partial_j$ omitted.

We assume that $f_1$ is of the form
$$
f_1=(a x_j+ b e_{1}^{(j)})\partial^{\alpha}\Delta,
$$
with $a$ and $b$ coefficients to be determined. Applying
$\nabla_1$ we get
$$
\nabla_1f_1=(a+(n-1)b)\partial^{\alpha}\Delta,
$$
while applying $\nabla_2$ we get
\begin{eqnarray*}
\nabla_2 f_1 & = & 2\sum_{i=1}^{n} \partial_i(a x_j+ b
e_{1}^{(j)})\partial_i\partial^{\alpha}\Delta\\
 & = & 2a\partial_j \partial^{\alpha}\Delta+2b
 \partial^{\alpha}\nabla_{1}^{(j)}\Delta\\
 & = & 2(a-b)\partial_j\partial^{\alpha}\Delta.
\end{eqnarray*}
Since the matrix
$$
\left(%
\begin{array}{cc}
  1 & n-1 \\
  1 & -1 \\
\end{array}%
\right)
$$
is invertible for every $n\geq 1$, we just showed how to construct
a solution of the system of equations
\begin{eqnarray*}
\nabla_1 f_1 & = & c\partial^{\alpha}\Delta,\\
\nabla_{2} f_1 & = & d\partial_j \partial^{\alpha}\Delta
\end{eqnarray*}
for any coefficients $c$ and $d$ and any $j$. All this together
with the observations in the first section takes care of the
$f_1$'s.

We indicate here how one could proceed to get an $f_2$ such that
\begin{eqnarray*}
\nabla_1 f_2 & = & -\widetilde{D}_1f_1\\
\nabla_{2} f_2 & = & -\widetilde{D}_2f_1,
\end{eqnarray*}
for $f_{1}^{(1)}=(ax_1+be_{1}^{(1)})\partial_{1}^2\Delta$. We have
\begin{eqnarray*}
-\widetilde{D}_1f_{1}^{(1)} & = &
-2\sum_{i=1}^{n}x_i\partial_i(ax_1+be_{1}^{(1)})\partial_i\partial_{1}^2\Delta\\
 & = & -2ax_1\partial_{1}^3\Delta
 -2b\sum_{i=2}^{n}x_i\partial_i \partial_{1}^2\Delta.
\end{eqnarray*}
Now if we set $g=\partial_{1}^{2}\prod_{i=2}^{n}(x_1-x_i)$ we have
\begin{eqnarray*}
\sum_{i=2}^{n}x_i\partial_i \partial_{1}^2\Delta & = & \left(\sum_{i=2}^{n}x_i\partial_i\Delta^{(1)}\right) g+\Delta^{(1)}\sum_{i=2}^{n}x_i\partial_ig\\
 & = & \left(\!\!\!%
\begin{array}{c}
  n-1 \\
  2 \\
\end{array}\!\!\!%
\right)\partial_{1}^2\Delta+(n-3)\partial_{1}^2\Delta+x_{1}\partial_{1}^3\Delta\\
 & = &
 \frac{n^2-n-4}{2}\partial_{1}^2\Delta+x_{1}\partial_{1}^3\Delta.
\end{eqnarray*}
Hence
$$
-\widetilde{D}_1f_{1}^{(1)}=-2(a+b)x_{1}\partial_{1}^3\Delta-b(n^2-n-4)\partial_{1}^2\Delta.
$$
Also
\begin{eqnarray*}
-\widetilde{D}_2f_{1}^{(1)} & = &
-3\sum_{i=1}^{n}x_i\partial_i(ax_1+be_{1}^{(1)})\partial_{i}^2\partial_{1}^2\Delta\\
 & = & -3ax_1\partial_{1}^4\Delta
 -3b\partial_{1}^2\left(\sum_{i=2}^{n}x_i\partial_{i}^2 \Delta\right).
\end{eqnarray*}
Now
\begin{eqnarray*}
\sum_{i=2}^{n}x_i\partial_{i}^2 \Delta & = &
2\sum_{i=2}^{n}\partial_i\Delta^{(1)} x_i\partial_ig\\
 & = & \Delta^{(1)}2P_1g =x_1\partial_{1}^2\Delta,
\end{eqnarray*}
hence
$$
-\widetilde{D}_2f_{1}^{(1)}=-3(a+b)x_1\partial_{1}^4\Delta-6b\partial_{1}^3\Delta.
$$
Since we already know how to take care of the terms
$-b(n^2-n-4)\partial_{1}^2\Delta$ and $-6b\partial_{1}^3\Delta$,
it will be more than enough to solve the following more general
problem:
\begin{eqnarray*}
\nabla_1 f_2 & = & \tilde{a}x_{1}\partial_{1}^k\Delta+\tilde{b} e_{1}^{(1)}\partial_{1}^k\Delta\\
\nabla_{2} f_2 & = &
\tilde{c}x_{1}\partial_{1}^{k+1}\Delta+\tilde{d}
e_{1}^{(1)}\partial_{1}^{k+1}\Delta,
\end{eqnarray*}
where $\hat{a},\hat{b},\hat{c}$ and $\hat{d}$ are coefficients,
and $k\geq 0$ is an integer.

Assume that $f_2$ is of the form
$$
f_2=(ax_{1}^2+bx_{1}e_{1}^{(1)}+c(e_{1}^{(1)})^2+de_{2}^{(1)})\partial_{1}^{k}\Delta+(\hat{a}x_1+\hat{b}e_{1}^{(1)})\partial_{1}^{k-1}\Delta,
$$
where $a,b,c,d,\hat{a}$ and $\hat{b}$ are coefficients to be
determined.

Now
$$
\nabla_1
f_2=((2a+(n-1)b)x_1+(b+2(n-1)c+(n-2)d)e_{1}^{(1)})\partial_{1}^{k}\Delta+(\hat{a}+(n-1)\hat{b})\partial_{1}^{k-1}\Delta,
$$
while
\begin{eqnarray*}
\nabla_2 f_2 & = &
(2a+2(n-1)c)\partial_{1}^{k}\Delta+2(\hat{a}-\hat{b})\partial_{1}^{k}\Delta\\
 & + & (4ax_1+2be_{1}^{(1)})\partial_{1}^{k+1}\Delta\\
 & + &
 (2bx_1+4ce_{1}^{(1)})\nabla_{1}^{(1)}\partial_{1}^{k}\Delta\\
 & + &
 d\partial_{1}^{k}\sum_{s=0}^{n-1}(-1)^{s+n-1}\nabla_2(e_{2}^{(1)}e_{n-1-s}^{(1)})x_{1}^s.
\end{eqnarray*}
Notice that the formula for $\nabla_2(e_ke_2)$ works also for
$k=1$. Hence the last term is
$$
d\left(2e_{1}^{(1)}\partial_{1}^{k+1}\Delta-2(n-1)\partial_{1}^{k}\Delta+2x_1\partial_{1}^{k+1}\Delta\right).
$$
Finally we have
\begin{eqnarray*}
\nabla_2 f_2 & = &
(2a+2(n-1)c-2(n-1)d+2(\hat{a}-\hat{b}))\partial_{1}^{k}\Delta\\
 & + & (4a-2b+2d)x_1\partial_{1}^{k+1}\Delta\\
 & + &
 (2b-4c+2d)e_{1}^{(1)}\partial_{1}^{k+1}\Delta.
\end{eqnarray*}
We already observed that with the coefficients $\hat{a}$ and
$\hat{b}$ we can get anything, hence we can disregard the terms
with $\partial_{1}^{k-1}\Delta$ and $\partial_{1}^{k}\Delta$.
What's left gives rise to a linear system with matrix
$$
\left(%
\begin{array}{cccc}
  2 & n-1 & 0 & 0 \\
  0 & 1 & 2(n-1) & (n-2) \\
  4 & -2 & 0 & 2 \\
  0 & 2 & -4 & 2 \\
\end{array}%
\right),
$$
whose determinant is $32(n^2-n)$. Hence for $n\geq 2$ this matrix
is non-singular, and this allows us to solve the system for all
values of $\widetilde{a},\widetilde{b},\widetilde{c}$ and
$\widetilde{d}$, and of course for any $k\geq 0$.

Using the Remark \ref{reductionrmk}, we can easily see that in
order to get an $f_2$ for any of the $f_1$ we found, we still need
to solve the system of equations
\begin{eqnarray*}
\nabla_1 f_2 & = & f_1\\
\nabla_{2} f_2 & = & \partial_j f_1.
\end{eqnarray*}
We have for $j\neq 1$
\begin{eqnarray*}
\partial_j f_1 & = & b\partial_{1}^2\Delta+(ax_1+be_{1}^{(1)})\partial_j\partial_{1}^2\Delta\\
 & = &
 b\partial_{1}^2\Delta-(ax_1+be_{1}^{(1)})\nabla_{1}^{(j)}\partial_{1}^2\Delta\\
 & = &
 b\partial_{1}^2\Delta-\nabla_{1}^{(j)}\left((ax_1+be_{1}^{(1)})\partial_{1}^2\Delta\right)+(a+(n-2)b)\partial_{1}^2\Delta.
\end{eqnarray*}
Also,
$$
\nabla_{1}^{(j)}\left((ax_1+be_{1}^{(1)})\partial_{1}^2\Delta\right)=\frac{1}{2}\nabla_{2}\left(e_{1}^{(j)}\cdot
(ax_1+be_{1}^{(1)})\partial_{1}^2\Delta
\right)-\frac{1}{2}e_{1}^{(j)}\partial_{1}^3\Delta.
$$
On the other hand,
$$
\nabla_1\left(e_{1}^{(j)}\cdot
(ax_1+be_{1}^{(1)})\partial_{1}^2\Delta
\right)=(n-1)(ax_1+be_{1}^{(1)})\partial_{1}^2\Delta+e_{1}^{(j)}\partial_{1}^2\Delta.
$$
Using what we have proved above, it's now clear that it's more
than enough to solve the system
\begin{eqnarray*}
\nabla_1 f_2 & = & \tilde{a}e_{1}^{(j)}\partial_{1}^{k}\Delta\\
\nabla_{2} f_2 & = & \tilde{b}e_{1}^{(j)}\partial_{1}^{k+1}\Delta,
\end{eqnarray*}
where $\tilde{a}$ and $\tilde{b}$ are arbitrary coefficients, and
$k\geq 0$ is an integer.

We leave the problem of finding a solution to this system open.

\section{Actions on alternating and symmetric polynomials}

We stick to the notation $e_k:=e_k(x_2,x_3,\dots,x_n)$, while
$e_{k}^{(i_1,i_2,\dots,i_r)}$ indicates the elementary symmetric
function of degree $k$ in the variables
$\{x_2,x_3,\dots,x_n\}\setminus \{i_1,i_2,\dots,i_r\}$. We recall
also the obvious relations
$$
\partial_je_{k}^{(i_1,i_2,\dots,i_r)}=e_{k-1}^{(i_1,i_2,\dots,i_r,j)},\quad
\text{and}\quad
e_{k}^{(i_1,i_2,\dots,i_r)}=e_{k}^{(i_1,i_2,\dots,i_r,j)}+x_je_{k-1}^{(i_1,i_2,\dots,i_r,j)}
$$
for $j\in \{x_2,x_3,\dots,x_n\}\setminus \{i_1,i_2,\dots,i_r\}$.

We remark also that all the identities that we are going to prove
will remain valid for elementary functions in any subset of the
variables involved, as long as we replace $n$ by the number of
variables involved plus one.

Another basic observation is that the elementary symmetric
functions $e_{\lambda}$'s and the $\Delta\cdot e_{\lambda}$'s,
where $\lambda$ runs over all partitions, form a basis of
symmetric and alternating polynomials respectively.

We are going to use all this without mentioning it anymore along
the way. Note also that we leave without proof the identities that
have been already proved in the previous sections.

In what follows $g$ will be a symmetric functions in the variables
$x_2,x_3,\dots,x_n$.

The action of $\nabla_1$ on symmetric functions is described by
the identity
$$
\nabla_1 e_k=(n-k)e_{k-1}
$$
together with Leibniz rule.

The action on alternating functions now follows immediately from
this one and Leibniz rule:
$$
\nabla_1(\Delta^{(1)} g)=(\nabla
\Delta^{(1)})g+\Delta^{(1)}(\nabla_1 g)=\Delta^{(1)}(\nabla_1 g).
$$

The following identity together with Leibniz rule describes the
action of the laplacian on symmetric functions.
\begin{Lem}
For $k\geq h$ we have
$$
\nabla_2(e_ke_h) =
2(n-k)e_{k-1}e_{h-1}-2\sum_{i=1}^{h-1}(k-h+2i)e_{k+i-1}e_{h-i-1}.
$$
\end{Lem}
\begin{proof}
We proceed by multiple induction on $k,h$ and $n$.
\begin{eqnarray*}
\nabla_2(e_ke_h) & = & 2\sum_{j=1}^{n}\partial_je_k
\partial_je_h\\
 & = &
 2\sum_{j=1}^{n}\partial_j(e_{k}^{(n)}+x_ne_{k-1}^{(n)})\cdot
 \partial_j(e_{h}^{(n)}+x_ne_{h-1}^{(n)})\\
 & = & 2\partial_n(e_{k}^{(n)}+x_ne_{k-1}^{(n)})\cdot
 \partial_n(e_{h}^{(n)}+x_ne_{h-1}^{(n)})\\
 & + &
 2\sum_{j=1}^{n-1}\partial_j(e_{k}^{(n)}+x_ne_{k-1}^{(n)})\cdot
 \partial_j(e_{h}^{(n)}+x_ne_{h-1}^{(n)})\qquad\qquad\qquad\qquad\\
 \quad & = & 2\,
 e_{k-1}^{(n)}e_{h-1}^{(n)}+2\sum_{j=1}^{n-1}\partial_je_{k}^{(n)}
\partial_je_{h}^{(n)}+2\sum_{j=1}^{n-1}\partial_jx_ne_{k}^{(n)}
\partial_je_{h-1}^{(n)}\\
 & + & 2\sum_{j=1}^{n-1}x_n\partial_je_{k-1}^{(n)}
\partial_je_{h}^{(n)}+2\sum_{j=1}^{n-1}x_{n}^2\partial_je_{k-1}^{(n)}
\partial_je_{h-1}^{(n)}\\
 & = & 2\,
 e_{k-1}^{(n)}e_{h-1}^{(n)}+2(n-k-1)e_{k-1}^{(n)}e_{h-1}^{(n)}
 -2\sum_{i=1}^{h-1}(k-h+2i)e_{k+i-1}^{(n)}e_{h-i-1}^{(n)}\\
 & + & x_n\left(2(n-k-1)e_{k-1}^{(n)}e_{h-2}^{(n)}
 -2\sum_{i=1}^{h-2}(k-h+1+2i)e_{k+i-1}^{(n)}e_{h-i-2}^{(n)}\right.\\
  & + & \left.2(n-k)e_{k-2}^{(n)}e_{h-1}^{(n)}
 -2\sum_{i=1}^{h-1}(k-h-1+2i)e_{k+i-2}^{(n)}e_{h-i-1}^{(n)}\right)\\
 & + & x_{n}^2\left(2(n-k)e_{k-2}^{(n)}e_{h-2}^{(n)}
 -2\sum_{i=1}^{h-2}(k-h+2i)e_{k+i-2}^{(n)}e_{h-i-2}^{(n)}\right)
\end{eqnarray*}
\begin{eqnarray*}
  & = & 2(n-k)e_{k-1}^{(n)}e_{h-1}^{(n)}
 -2\sum_{i=1}^{h-1}(k-h+2i)e_{k+i-1}^{(n)}e_{h-i-1}^{(n)}\\
 & + &
 x_n\left(2(n-k)\left(e_{k-1}^{(n)}e_{h-2}^{(n)}+e_{k-2}^{(n)}e_{h-1}^{(n)}\right)
 \right.\\
  & - &
  \left.2\sum_{i=1}^{h-1}(k-h+2i)\left(e_{k+i-1}^{(n)}e_{h-i-2}^{(n)}+e_{k+i-2}^{(n)}e_{h-i-1}^{(n)}\right)\right)\\
 & + & x_{n}^2\left(2(n-k)e_{k-2}^{(n)}e_{h-2}^{(n)}
 -2\sum_{i=1}^{h-2}(k-h+2i)e_{k+i-2}^{(n)}e_{h-i-2}^{(n)}\right)\\
 & = &
 2(n-k)e_{k-1}e_{h-1}-2\sum_{i=1}^{h-1}(k-h+2i)e_{k+i-1}e_{h-i-1}.
\end{eqnarray*}
The base cases are trivial.
\end{proof}

The action of the laplacian on alternating functions now follows
from
$$
\frac{1}{\Delta^{(1)}} \nabla_2(\Delta^{(1)} g)=(\nabla_2+2P_2) g,
$$
where
$$
P_2:=\sum_{2\leq i<j\leq
n}\frac{1}{x_i-x_j}(\partial_i-\partial_j),
$$
the formula
$$
P_2 e_k  =  -\left(\!\!\!\begin{array}{c}
 n-k+1 \\
  2 \\
\end{array} \!\!\!\right) e_{k-2},
$$
and Leibniz rule.

The following identity together with Leibniz rule describes the
action of the operator $\widetilde{D}_1$ on symmetric functions.
\begin{Lem}
For $k\geq h$,
$$
\widetilde{D}_1(e_ke_h)=2\sum_{i=0}^{h-1}(k-h+1+2i)e_{k+i}e_{h-1-i}.
$$
\end{Lem}
\begin{proof}
We proceed by multiple induction on $k,h$ and $n$.
\begin{eqnarray*}
\widetilde{D}_1(e_ke_h) & = & 2\sum_{i=2}^{n}x_i\partial_i e_k
\partial_i e_h\\
 & = & 2\sum_{i=2}^{n}x_i\partial_i
 (e_{k}^{(n)}+x_ne_{k-1}^{(n)})\partial_i
 (e_{h}^{(n)}+x_ne_{h-1}^{(n)})\\
 & = & 2\, x_ne_{k-1}^{(n)}e_{h-1}^{(n)}+2\sum_{i=2}^{n-1}x_i\partial_i
 (e_{k}^{(n)}+x_ne_{k-1}^{(n)})\partial_i
 (e_{h}^{(n)}+x_ne_{h-1}^{(n)})\\
 & = & 2\, x_ne_{k-1}^{(n)}e_{h-1}^{(n)}+2\sum_{i=2}^{n-1}x_i\partial_i
 e_{k}^{(n)}
\partial_i e_{h}^{(n)}\\
 & + & 2\, x_n\sum_{i=2}^{n-1}x_i\left(\partial_i
 e_{k}^{(n)}
\partial_i e_{h-1}^{(n)}+\partial_i
 e_{k-1}^{(n)}
\partial_i e_{h}^{(n)}\right)\\
 & + & 2\,x_{n}^2\sum_{i=2}^{n-1}x_i\partial_i
 e_{k-1}^{(n)}
\partial_i e_{h-1}^{(n)}\\
 & = & 2\, x_ne_{k-1}^{(n)}e_{h-1}^{(n)}+2\sum_{i=0}^{h-1}(k-h+1+2i)e_{k+i}^{(n)}e_{h-1-i}^{(n)}\\
 & + & 2\, x_n\left(\sum_{i=0}^{h-2}(k-h+2+2i)e_{k+i}^{(n)}e_{h-2-i}^{(n)}+\sum_{i=0}^{h-1}(k-h+2i)e_{k-1+i}^{(n)}e_{h-1-i}^{(n)}\right)\\
 & + &
 2\,x_{n}^2\sum_{i=0}^{h-2}(k-h+1+2i)e_{k-1+i}^{(n)}e_{h-2-i}^{(n)}
\end{eqnarray*}
 \begin{eqnarray*}
  & = & 2\sum_{i=0}^{h-1}(k-h+1+2i)e_{k+i}^{(n)}e_{h-1-i}^{(n)}\\
 & + & 2\, x_n\left(\sum_{i=0}^{h-2}(k-h+1+2i)e_{k+i}^{(n)}e_{h-2-i}^{(n)}+\sum_{i=0}^{h-1}(k-h+1+2i)e_{k-1+i}^{(n)}e_{h-1-i}^{(n)}\right)\\
 & + &
 2\,x_{n}^2\sum_{i=0}^{h-2}(k-h+1+2i)e_{k-1+i}^{(n)}e_{h-2-i}^{(n)}\\
 & = & 2\sum_{i=0}^{h-1}(k-h+1+2i)e_{k+i}e_{h-1-i}.
\end{eqnarray*}
The base cases are trivial.
\end{proof}

We have
$$
\frac{1}{\Delta^{(1)}} \widetilde{D}_1(\Delta^{(1)}
g)=(2P_1+\widetilde{D}_1)g,
$$
where
$$
P_1:=\sum_{2\leq i<j\leq
n}\frac{1}{x_i-x_j}(x_i\partial_i-x_j\partial_j).
$$
We have the following identity, whose proof is analogous to the
one of the identities (\ref{p2}):
$$
P_1 e_k  =  \left(\!\!\!\begin{array}{c}
 n-k \\
  2 \\
\end{array}\!\!\! \right) e_{k-1}.
$$
All this together with Leibniz rule describes the action of
$\widetilde{D}_1$ on alternating polynomials.

The following identity together with Leibniz rule describes the
action of the operator $\widetilde{D}_2$ on symmetric functions.
\begin{Lem}
For $k\geq h\geq l$,
$$
\widetilde{D}_2(e_ke_he_l)=6\left(\sum_{j=0}^{l-1}\sum_{i=0}^{h-1}(k-h+1+j+2i)e_{k+i+j}e_{h-1-i}e_{l-1-j}\right.\qquad\qquad\qquad
$$
$$
\qquad\qquad\qquad\qquad\qquad\qquad\qquad
-\left.\sum_{j=0}^{l-2}\sum_{i=1}^{l-1-j}(h-l+j+2i)e_{k+j}e_{h-1+i}e_{l-1-i-j}\right).
$$
\end{Lem}
\begin{proof}
We proceed by multiple induction on $k,h,l$ and $n$.
\begin{eqnarray*}
\frac{1}{6}\widetilde{D}_2(e_ke_he_l) & = &
\sum_{i=2}^{n}x_i\partial_ie_k \partial_ie_h \partial_ie_l\\
 & = & \sum_{i=2}^{n}x_i\partial_i(e_{k}^{(n)}+x_n e_{k-1}^{(n)}) \partial_i(e_{h}^{(n)}+x_n e_{h-1}^{(n)}) \partial_i(e_{l}^{(n)}+x_n
 e_{l-1}^{(n)})
\end{eqnarray*}
\begin{eqnarray*}
 & = & x_n e_{k-1}^{(n)} e_{h-1}^{(n)} e_{l-1}^{(n)} +
 \sum_{i=2}^{n-1}x_i\partial_ie_{k}^{(n)} \partial_ie_{h}^{(n)}
 \partial_ie_{l}^{(n)}\\
 & + & x_n\left(\sum_{i=2}^{n-1}x_i\left(\partial_ie_{k-1}^{(n)} \partial_ie_{h}^{(n)}
 \partial_ie_{l}^{(n)}+\partial_ie_{k}^{(n)} \partial_ie_{h-1}^{(n)}
 \partial_ie_{l}^{(n)}+\partial_ie_{k}^{(n)} \partial_ie_{h}^{(n)}
 \partial_ie_{l-1}^{(n)}  \right)\right)\\
 & + & x_{n}^2\left(\sum_{i=2}^{n-1}x_i\left(\partial_ie_{k}^{(n)} \partial_ie_{h-1}^{(n)}
 \partial_ie_{l-1}^{(n)}+\partial_ie_{k-1}^{(n)} \partial_ie_{h}^{(n)}
 \partial_ie_{l-1}^{(n)}+\partial_ie_{k-1}^{(n)} \partial_ie_{h-1}^{(n)}
 \partial_ie_{l}^{(n)}  \right)\right)\\
 & + & x_{n}^3\sum_{i=2}^{n-1}x_i\partial_ie_{k-1}^{(n)} \partial_ie_{h-1}^{(n)}
 \partial_ie_{l-1}^{(n)}.
\end{eqnarray*}
At this point we use induction, replacing the suitable terms by
our formula. To be more efficient, we analyze the expansion with
respect to powers of $x_n$.

For the factor of $x_n$ we get
\begin{eqnarray*}
e_{k-1}^{(n)} e_{h-1}^{(n)} e_{l-1}^{(n)} & + &
\sum_{j=0}^{l-1}\sum_{i=0}^{h-1}
 (k-h+j+2i)e_{k-1+i+j}^{(n)}e_{h-1-i}^{(n)}e_{l-1-j}^{(n)}\\
 & + & \sum_{j=0}^{l-1}\sum_{i=0}^{h-2}
 (k-h+2+j+2i)e_{k+i+j}^{(n)}e_{h-2-i}^{(n)}e_{l-1-j}^{(n)}\\
& + &\sum_{j=0}^{l-2}\sum_{i=0}^{h-1}
(k-h+1+j+2i)e_{k+i+j}^{(n)}e_{h-1-i}^{(n)}e_{l-2-j}^{(n)}\\
& - & \sum_{j=0}^{l-2}\sum_{i=1}^{l-1-j}
(h-l+j+2i)e_{k-1+j}^{(n)}e_{h-1+i}^{(n)}e_{l-1-i-j}^{(n)}\\
& - & \sum_{j=0}^{l-2}\sum_{i=1}^{l-1-j}
(h-l-1+j+2i)e_{k+j}^{(n)}e_{h-2+i}^{(n)}e_{l-1-i-j}^{(n)}\\
& - & \sum_{j=0}^{l-3}\sum_{i=1}^{l-2-j}
(h-l+1+j+2i)e_{k+j}^{(n)}e_{h-1+i}^{(n)}e_{l-2-i-j}^{(n)}.
\end{eqnarray*}
Rearranging the terms we get what we want:
$$
\sum_{j=0}^{l-2}\sum_{i=0}^{h-1}
(k-h+1+j+2i)\left(e_{k+i+j}^{(n)}e_{h-2-i}^{(n)}e_{l-1-j}^{(n)}+e_{k+i+j}^{(n)}e_{h-1-i}^{(n)}e_{l-2-j}^{(n)}
+e_{k-1+i+j}^{(n)}e_{h-1-i}^{(n)}e_{l-1-j}^{(n)}\right)+
$$
$$
- \sum_{j=0}^{l-2}\sum_{i=1}^{l-1-j}
(h-l+j+2i)\left(e_{k-1+j}^{(n)}e_{h-1+i}^{(n)}e_{l-1-i-j}^{(n)}+e_{k+j}^{(n)}e_{h-2+i}^{(n)}e_{l-1-i-j}^{(n)}+
e_{k+j}^{(n)}e_{h-1+i}^{(n)}e_{l-2-i-j}^{(n)}\right)
$$

Analogously for the factor of $x_{n}^2$. What is left is already
what we want.

The base cases are trivial.
\end{proof}

We have
$$
\frac{1}{\Delta^{(1)}} \widetilde{D}_2(\Delta^{(1)}
g)=(6Q_2+3\widetilde{P}_2+\widetilde{D}_2)g,
$$
where
$$
Q_2:=\sum_{j=1}^n\sum_{i<k}^{(j)}\frac{1}{(x_j-x_i)(x_j-x_k)}x_j\partial_{j},
$$
and
$$
\widetilde{P}_2:=\sum_{1\leq i<j\leq
n}\frac{1}{x_i-x_j}(x_i\partial_{i}^2-x_j\partial_{j}^2).
$$
The following Lemma together with Leibniz rule describes the
action of $Q_2$ on symmetric polynomials.
\begin{Lem}
We have
$$
Q_2 e_k  =  -\left(\!\!\!\begin{array}{c}
 n-k+1 \\
  3 \\
\end{array}\!\!\! \right) e_{k-2}.
$$
\end{Lem}
\begin{proof}
It's clear that we have the following relations:
$$
e_{m}^{(i_1,i_2,\dots,i_r)}=e_{m}^{(i_1,i_2,\dots,i_r,j)}+x_je_{m-1}^{(i_1,i_2,\dots,i_r,j)},
$$
for all $j\notin \{i_1,\dots,i_r\}$. We are going to use them
repeatedly without mentioning it.

For $2 \leq i<j<k\leq n$ we have
$$
\frac{x_j\partial_je_{m}}{(x_j-x_i)(x_j-x_k)}+\frac{x_i\partial_ie_{m}}{(x_i-x_j)(x_i-x_k)}+
\frac{x_k\partial_ke_{m}}{(x_k-x_i)(x_k-x_j)}=\qquad\qquad\qquad\qquad\qquad
$$
\begin{eqnarray*}
 \qquad & = &
 \frac{-x_j(x_i-x_k)e_{m-1}^{(j)}+x_i(x_j-x_k)e_{m-1}^{(i)}+x_k(x_i-x_j)e_{m-1}^{(k)}}{(x_i-x_j)(x_j-x_k)(x_i-x_k)}\\
  & = &
  -\left(\frac{x_j(x_k-x_i)e_{m-1}^{(j)}+x_i(x_j-x_k)e_{m-1}^{(i)}+x_k(x_i-x_j)e_{m-1}^{(k)}}{(x_i-x_j)(x_j-x_k)(x_k-x_i)}\right).
\end{eqnarray*}
Clearly the denominator divides the numerator, but we want to
compute the quotient. The numerator is equal to
$$
(x_k-x_i)(x_je_{m-1}^{(j)})+x_ix_k(e_{m-1}^{(k)}-e_{m-1}^{(i)})+x_j(x_ie_{m-1}^{(i)}-x_ke_{m-1}^{(k)})=\qquad\qquad\qquad
$$
\begin{eqnarray*}
\qquad\qquad & = &
(x_k-x_i)(x_je_{m-1}^{(j)})+x_ix_k(x_i-x_k)e_{m-2}^{(i,k)}+x_j(x_i-x_k)e_{m-1}^{(i,k)}\\
 & = &
 (x_k-x_i)(x_je_{m-1}^{(j)}-x_ix_ke_{m-2}^{(i,k)}-x_je_{m-1}^{(i,k)}).
\end{eqnarray*}
The second factor of the last term is equal to
$$
x_je_{m-1}^{(i,j,k)}+x_jx_ie_{m-2}^{(i,j,k)}+x_jx_ke_{m-2}^{(i,j,k)}
+ x_ix_jx_k e_{m-3}^{(i,j,k)}+\qquad\qquad\qquad
$$
\begin{eqnarray*}
\qquad\qquad\qquad
 & - &
 x_ix_ke_{m-2}^{(i,j,k)}-x_ix_jx_ke_{m-3}^{(i,j,k)}-x_je_{m-1}^{(i,j,k)}-x_{j}^2e_{m-2}^{(i,j,k)}\\
  & = & (x_ix_j+x_jx_k-x_ix_k-x_{j}^2)e_{m-2}^{(i,j,k)}\\
  & = & (x_i-x_j)(x_j-x_k)e_{m-2}^{(i,j,k)}.
\end{eqnarray*}
In conclusion we get
$$
Q_2 e_m= - \sum_{2\leq i<j<k\leq
n}e_{m-2}^{(i,j,k)}=-\left(\!\!\!\begin{array}{c}
 n-m+1 \\
  3 \\
\end{array}\!\!\! \right)e_{m-2},
$$
where the last equality comes from counting how many times the
monomial $x_2x_3\cdots x_{m-1}$ shows up.
\end{proof}

The following identity together with Leibniz rule describes the
action of the operator $\widetilde{P}_2$ on symmetric functions.
\begin{Lem}
For $k\geq h$,
$$
\widetilde{P}_2(e_ke_h)=(n-k)(n-k-1)e_{k-1}e_{h-1}-(2n-h-k-1)\left(\sum_{i=1}^{h-1}(k-h+2i)e_{k-1+i}e_{h-1-i}\right).
$$
\end{Lem}
\begin{proof}
By induction on $n$:
\begin{eqnarray*}
\widetilde{P}_2(e_ke_h)\!\!\!\!\! & =\!\! & 2\sum_{2\leq i<j\leq
n}\frac{1}{x_i-x_j}\left(x_i\partial_i e_k\partial_i
e_h-x_j\partial_j  e_k\partial_j e_h\right)\\
 & =\!\! & 2\sum_{2\leq i<j\leq n}\frac{1}{x_i-x_j}\left(x_i
e_{k-1}^{(i)}e_{h-1}^{(i)}-x_j e_{k-1}^{(j)}
e_{h-1}^{(j)}\right)\\
 & =\!\! & 2\sum_{2\leq i<j\leq n}\frac{1}{x_i-x_j}\left(x_i
\left(e_{k-1}^{(i,j)}+x_je_{k-2}^{(i,j)}\right)\left(e_{h-1}^{(i,j)}+x_je_{h-2}^{(i,j)}\right)\right.\\
 & -\!\! & \left.x_j
\left(e_{k-1}^{(i,j)}+x_ie_{k-2}^{(i,j)}\right)\left(e_{h-1}^{(i,j)}+x_ie_{h-2}^{(i,j)}\right)\right)\\
 & =\!\! & 2\sum_{2\leq i<j\leq n}\frac{1}{x_i-x_j}\left(
 (x_i-x_j)e_{k-1}^{(i,j)}e_{h-1}^{(i,j)}+x_ix_j(x_j-x_i)e_{k-2}^{(i,j)}e_{h-2}^{(i,j)}\right)\\
 & =\!\! & 2\sum_{2\leq i<j\leq n}\left(
 e_{k-1}^{(i,j)}e_{h-1}^{(i,j)}-x_ix_je_{k-2}^{(i,j)}e_{h-2}^{(i,j)}\right)\\
 & =\!\! & 2\sum_{2\leq i<n}\left(
 e_{k-1}^{(i,n)}e_{h-1}^{(i,n)}-x_ix_ne_{k-2}^{(i,n)}e_{h-2}^{(i,n)}\right)\\
 & +\!\! & 2\sum_{2\leq i<j\leq n-1}\left(
 e_{k-1}^{(i,j,n)}e_{h-1}^{(i,j,n)}-x_ix_je_{k-2}^{(i,j,n)}e_{h-2}^{(i,j,n)}\right)\\
  & +\!\! & 2\, x_n\sum_{2\leq i<j\leq n-1}\left(\left(
 e_{k-1}^{(i,j,n)}e_{h-2}^{(i,j,n)}+e_{k-2}^{(i,j,n)}e_{h-1}^{(i,j,n)}\right)-x_ix_j\left(e_{k-2}^{(i,j,n)}e_{h-3}^{(i,j,n)}
 +e_{k-3}^{(i,j,n)}e_{h-2}^{(i,j,n)}\right)\right)\\
 & +\!\! & 2\, x_{n}^2\sum_{2\leq i<j\leq n-1}\left(
 e_{k-2}^{(i,j,n)}e_{h-2}^{(i,j,n)}-x_ix_je_{k-3}^{(i,j,n)}e_{h-3}^{(i,j,n)}\right)
 \end{eqnarray*}
 \begin{eqnarray*}
  & =\!\! & 2\sum_{2\leq i<n}
 \partial_i e_{k}^{(n)}\partial_ie_{h}^{(n)}- 2\, x_n\sum_{2\leq i<n}x_i \partial_i e_{k-1}^{(n)}\partial_i e_{h-1}^{(n)}\\
 & +\!\! & (n-k-1)(n-k-2)e_{k-1}^{(n)}e_{h-1}^{(n)}-(2n-h-k-3)\left(\sum_{i=1}^{h-1}(k-h+2i)e_{k-1+i}^{(n)}e_{h-1-i}^{(n)}\right)\\
 & +\!\! & x_n\left(\!\!(n-k-1)(n-k-2)e_{k-1}^{(n)}e_{h-2}^{(n)}-(2n-h-k-2)\left(\sum_{i=1}^{h-2}(k-h+1+2i)e_{k-1+i}^{(n)}e_{h-2-i}^{(n)}\right)
 \!\!\right)\\
 & + \!\!& x_n\left(\!\!(n-k)(n-k-1)e_{k-2}^{(n)}e_{h-1}^{(n)}-(2n-h-k-2)\left(\sum_{i=1}^{h-1}(k-h-1+2i)e_{k-2+i}^{(n)}e_{h-1-i}^{(n)}\right)
 \!\!\right)\\
 & + \!\!&
 x_{n}^2\left(\!\!(n-k)(n-k-1)e_{k-2}^{(n)}e_{h-2}^{(n)}-(2n-h-k-1)\left(\sum_{i=1}^{h-2}(k-h+2i)e_{k-2+i}^{(n)}e_{h-2-i}^{(n)}\right)\!\!\right).
\end{eqnarray*}
We have
$$
2\sum_{2\leq i<n}
 \partial_i e_{k}^{(n)}\partial_ie_{h}^{(n)}- 2\, x_n\sum_{2\leq i<n}x_i \partial_i e_{k-1}^{(n)}\partial_i
 e_{h-1}^{(n)}= \nabla_2 (e_{k}^{(n)}e_{h}^{(n)})-
 x_n\widetilde{D}_1(e_{k-1}^{(n)}e_{h-1}^{(n)}),
$$
Hence
\begin{eqnarray*}
 \widetilde{P}_2(e_ke_h)\!\!\!\! & = \!\!& 2(n-k-1)e_{k-1}^{(n)}e_{h-1}^{(n)}-2\sum_{i=1}^{h-1}(k-h+2i)e_{k+i-1}^{(n)}e_{h-i-1}^{(n)} \\
   & - \!\!&   x_n\left(2\sum_{i=0}^{h-2}(k-h+1+2i)e_{k-1+i}^{(n)}e_{h-2-i}^{(n)}\right)\\
 & +\!\! & (n-k-1)(n-k-2)e_{k-1}^{(n)}e_{h-1}^{(n)}-(2n-h-k-3)\left(\sum_{i=1}^{h-1}(k-h+2i)e_{k-1+i}^{(n)}e_{h-1-i}^{(n)}\right)\\
 & +\!\! & x_n\left(\!\!(n-k-1)(n-k-2)e_{k-1}^{(n)}e_{h-2}^{(n)}-(2n-h-k-2)\left(\sum_{i=1}^{h-2}(k-h+1+2i)e_{k-1+i}^{(n)}e_{h-2-i}^{(n)}\right)
 \!\!\right)\\
 & + \!\!& x_n\left(\!\!(n-k)(n-k-1)e_{k-2}^{(n)}e_{h-1}^{(n)}-(2n-h-k-2)\left(\sum_{i=1}^{h-1}(k-h-1+2i)e_{k-2+i}^{(n)}e_{h-1-i}^{(n)}\right)
 \!\!\right)\\
 & +\!\! & x_{n}^2\left(\!\!(n-k)(n-k-1)e_{k-2}^{(n)}e_{h-2}^{(n)}-(2n-h-k-1)\left(\sum_{i=1}^{h-2}(k-h+2i)e_{k-2+i}^{(n)}e_{h-2-i}^{(n)}\right)
 \!\!\right)
\end{eqnarray*}
\begin{eqnarray*}
 \qquad\quad & = \!\!& (n-k)(n-k-1)e_{k-1}^{(n)}e_{h-1}^{(n)}-(2n-h-k-1)\left(\sum_{i=1}^{h-1}(k-h+2i)e_{k-1+i}^{(n)}e_{h-1-i}^{(n)}\right)\\
 & +\!\! & x_n\left(\!\!(n-k)(n-k-1)e_{k-1}^{(n)}e_{h-2}^{(n)}-(2n-h-k-1)\left(\sum_{i=1}^{h-2}(k-h+1+2i)e_{k-1+i}^{(n)}e_{h-2-i}^{(n)}\right)
 \!\!\right)\\
 & + \!\!& x_n\left(\!\!(n-k)(n-k-1)e_{k-2}^{(n)}e_{h-1}^{(n)}-(2n-h-k-1)\left(\sum_{i=1}^{h-1}(k-h-1+2i)e_{k-2+i}^{(n)}e_{h-1-i}^{(n)}\right)
 \!\!\right.\\
 & + \!\!& \left.((2n-h-k-1)-2(n-k-1)-(k-h+1))e_{k-1}^{(n)}e_{h-2}^{(n)} \right)\\
 & +\!\! & x_{n}^2\left(\!\!(n-k)(n-k-1)e_{k-2}^{(n)}e_{h-2}^{(n)}-(2n-h-k-1)\left(\sum_{i=1}^{h-2}(k-h+2i)e_{k-2+i}^{(n)}e_{h-2-i}^{(n)}\right)
 \!\!\right)\\
 & = &
 (n-k)(n-k-1)e_{k-1}e_{h-1}-(2n-h-k-1)\left(\sum_{i=1}^{h-1}(k-h+2i)e_{k-1+i}e_{h-1-i}\right).
\end{eqnarray*}
The base cases are trivial.
\end{proof}
All this together with Leibniz rule describes the action of
$\widetilde{D}_2$ on alternating polynomials.

\subsection{List of Formulae}

For convenience and for future reference, we give a list of the
formulae that we found along the way. \textsl{\textbf{In this
subsection}} we state them in terms of the variables
$x_1,x_2,\dots,x_n$, adapting the definitions accordingly.\newline

\textit{Here $e_k$ will be the elementary symmetric function in
$n$ variables of degree $k$, and $g$ a symmetric function in the
variables $x_1,x_2,\dots,x_n$.}

\subsubsection{Action of $\nabla_1$}
$$
\frac{1}{\Delta} \nabla_1(\Delta g)=\nabla_1 g.
$$
$$
\nabla_1 e_k=(n-k+1)e_{k-1}.
$$

\subsubsection{Action of $\nabla_2$}
$$
\frac{1}{\Delta} \nabla_2(\Delta g)=(\nabla_2+2P_2) g,
$$
where
$$
P_2:=\sum_{1\leq i<j\leq
n}\frac{1}{x_i-x_j}(\partial_i-\partial_j).
$$
$$
P_2 e_k  =  -\left(\!\!\!\begin{array}{c}
 n-k+2 \\
  2 \\
\end{array} \!\!\!\right) e_{k-2}.
$$
For $k\geq h$ we have
$$
\nabla_2(e_ke_h) =
2(n-k+1)e_{k-1}e_{h-1}-2\sum_{i=1}^{h-1}(k-h+2i)e_{k+i-1}e_{h-i-1}.
$$
If $g$ is a symmetric function,
$$
[\nabla_1,P_2]g=0.
$$

\subsubsection{Action of $\widetilde{D}_1$}

$$
\frac{1}{\Delta} \widetilde{D}_1(\Delta
g)=(2P_1+\widetilde{D}_1)g,
$$
where
$$
P_1:=\sum_{1\leq i<j\leq
n}\frac{1}{x_i-x_j}(x_i\partial_i-x_j\partial_j).
$$
$$
P_1 e_k  =  \left(\!\!\!\begin{array}{c}
 n-k+1 \\
  2 \\
\end{array}\!\!\! \right) e_{k-1}.
$$
For $k\geq h$,
$$
\widetilde{D}_1(e_ke_h)=2\sum_{i=0}^{h-1}(k-h+1+2i)e_{k+i}e_{h-1-i}.
$$

\subsubsection{Action of $\widetilde{D}_2$}

$$
\frac{1}{\Delta} \widetilde{D}_2(\Delta
g)=(6Q_2+3\widetilde{P}_2+\widetilde{D}_2)g,
$$
where
$$
Q_2:=\sum_{j=1}^n\sum_{i<k}^{(j)}\frac{1}{(x_j-x_i)(x_j-x_k)}x_j\partial_{j},
$$
and
$$
\widetilde{P}_2:=\sum_{1\leq i<j\leq
n}\frac{1}{x_i-x_j}(x_i\partial_{i}^2-x_j\partial_{j}^2).
$$
$$
Q_2 e_k  =  -\left(\!\!\!\begin{array}{c}
 n-k+2 \\
  3 \\
\end{array}\!\!\! \right) e_{k-2}.
$$
For $k\geq h$,
$$
\widetilde{P}_2(e_ke_h)=(n-k+1)(n-k)e_{k-1}e_{h-1}-(2n-h-k+1)\left(\sum_{i=1}^{h-1}(k-h+2i)e_{k-1+i}e_{h-1-i}\right).
$$
For $k\geq h\geq l$,
$$
\widetilde{D}_2(e_ke_he_l)=6\left(\sum_{j=0}^{l-1}\sum_{i=0}^{h-1}(k-h+1+j+2i)e_{k+i+j}e_{h-1-i}e_{l-1-j}\right.\qquad\qquad\qquad
$$
$$
\qquad\qquad\qquad\qquad\qquad\qquad\qquad
-\left.\sum_{j=0}^{l-2}\sum_{i=1}^{l-1-j}(h-l+j+2i)e_{k+j}e_{h-1+i}e_{l-1-i-j}\right).
$$

\section{Singular $q_0$-harmonics}

\textit{\textbf{Warning}: in this section we use the notations of
section 5.1}\newline

Recall here that in \cite{hivert} Thi\'ery and Hivert stated the
conjecture that in the case where $q$ is a complex number not of
the form $-a/b$ where $a\in \{1,2,\dots ,n\}$ and $b\in
\mathbb{N}$, we have the equality
$$
\sum_{d\geq 0}\dim  \pi_d(\mathcal{H}_{\mathbf{x};q})t^d = [n]_t
!.
$$
Inspired by a similar definition in \cite{hivert}, we define a
complex number $q_0$ \textit{singular} if the Frobenius
characteristic $F_{n;q_0}(t)$ of the $q_0$-harmonics is different
from the Frobenius characteristic $F_{n}(t)=F_{n;0}(t)$ of the
classical harmonics, which is (see \cite{macdonald})
$$
F_{n}(t)=\sum_{\lambda \vdash n}s_{\lambda}\sum_{T\in ST(n)
}t^{co(T)},
$$
where $\lambda \vdash n$ indicates that $\lambda$ is a partition
of $n$, $s_{\lambda}$ is the Schur function indexed by $\lambda$,
$ST(\lambda)$ denotes the set of standard tableaux of shape
$\lambda$, and $co(T)$ denote the cocharge of the tableau $T$.

One of the main result of this section is the following theorem.
\begin{Thm} \label{thmmain}
The values of $q_0$ of the form $-a/b$ where $a\in \{1,2,\dots
,n\}$, $b\in \mathbb{N}$ and $b\geq n$ are singular.
\end{Thm}
\begin{Remrk}
Notice that in the statement we don't require that $a$ and $b$ are
coprime. For example if $n=6$, then we will show that $-2/3$ is
singular, since it can be written as $-4/6$.
\end{Remrk}
More generally, in this appendix we will investigate the
$q_0$-harmonics for singular values of $q_0$.

\begin{Remrk}
Since the case $q_0=0$ reduces to the well known case of classical
$\frak{S}_n$-harmonics, in this section we will always assume
$q_0\neq 0$. Recall also, from the easy relations
$$
[D_{k;q_0},D_{h;q_0}]=q_0(k-h)D_{k+h;q_0},
$$
it follows that a polynomial $f$ is in
$\mathcal{H}_{\mathbf{x};q_0}$ if and only if
$$
D_{1;q_0}f=D_{2;q_0}f=0.
$$
We will use repeatedly this observation without mentioning it
anymore.
\end{Remrk}

In our computer investigations we realized that polynomials of
certain forms are $q_0$-harmonics for special values of $q_0$.
Using the formulae of the previous section we are now able to
prove that this is the case.

First of all we prove that for $1\leq k<n$ and $q_0=-1/(n-k)$ the
alternant $\Delta e_k$ is in $\mathcal{H}_{\mathbf{x};q_0}$. This
shows immediately that these values of $q_0$ are singular, since
in the classical case the only alternant is $\Delta$ in degree
$\left(\!\!\!\begin{array}{c}
 n \\
  2 \\
\end{array}\!\!\! \right)$.
\begin{Thm} \label{thmalter}
The polynomial $\Delta e_{k}$ is $q_0$-harmonic if and only if
$k<n$ and $q_0=-1/(n-k)$.
\end{Thm}
\begin{proof}
Let's look at the action of $D_{1;q_0}=\nabla_1+q_0
\widetilde{D}_1$ on $\Delta e_{k}$. Using the formulae listed in
the previous section, we have
\begin{eqnarray*}
D_{1;q_0} \Delta e_k & = & (\nabla_1+q_0\widetilde{D}_1) \Delta
e_k=\Delta (\nabla_1+q_0(2P_1+\widetilde{D}_1))
e_k\\
& = & \left((n-k+1)+2q_0\left(\!\!\!\begin{array}{c}
 n-k+1 \\
  2 \\
\end{array}\!\!\! \right)\right)\Delta e_k.
\end{eqnarray*}
Hence to have $D_{1;q_0}\Delta e_k=0$ we need to have $k<n$ and
$$
q_0=-\frac{1}{n-k}.
$$
Let's now look at $D_{2;q_0}\Delta e_k$. We have
\begin{eqnarray*}
D_{2;q_0}\Delta e_k & = &
\Delta((\nabla_2+2P_2)+q_0(6Q_2+3\widetilde{P}_2+\widetilde{D}_2))e_k\\
 & = & \left(-2\left(\!\!\!\begin{array}{c}
 n-k+2 \\
  2 \\
\end{array}\!\!\! \right)-q_0 6\left(\!\!\!\begin{array}{c}
 n-k+2 \\
  3 \\
\end{array}\!\!\! \right)\right) \Delta e_k,
\end{eqnarray*}
which is $0$ for $q_0=-1/(n-k)$.
\end{proof}
We determine now another class of $q_0$-harmonics which will imply
the singularity of many values of $q_0$.

Recall that we work in $n\geq 2$ variables.

\begin{Thm}
The polynomial $e_{1}^{m}(x_1,x_2,\dots,x_k)(x_1-x_2)$, with
$2\leq k\leq n$ and $m\geq 1$ is a $q_0$-harmonic if and only if
$q_0=-\frac{k}{m+1}$.
\end{Thm}
\begin{proof}
Let's look at the action of $D_{1;q_0}=\nabla_1+q_0
\widetilde{D}_1$. We have
$$
\nabla_1e_{1}^{m}(x_1,x_2,\dots,x_k)(x_1-x_2)=m\,k\,
e_{1}^{m-1}(x_1,x_2,\dots,x_k)(x_1-x_2),
$$
while
$$
\widetilde{D}_1e_{1}^{m}(x_1,x_2,\dots,x_k)(x_1-x_2) =
\qquad\qquad\qquad\qquad\qquad $$

\begin{eqnarray*}
\qquad\qquad\qquad & = & \left(2 \left(\!\!\!\begin{array}{c}
 m \\
  2 \\
\end{array}\!\!\!
\right)+2m\right)e_{1}^{m-1}(x_1,x_2,\dots,x_k)(x_1-x_2)\\
 & = &  m(m+1)\, e_{1}^{m-1}(x_1,x_2,\dots,x_k)(x_1-x_2).
\end{eqnarray*}
Therefore
$$
D_{1;q_0}e_{1}^{m}(x_1,x_2,\dots,x_k)(x_1-x_2)=\left(k\, m +q_0
m(m+1)\right)e_{1}^{m-1}(x_1,x_2,\dots,x_k)(x_1-x_2),
$$
and this is equal to $0$ if and only if $q_0=-\frac{k}{m+1}$.

We are left to check that also $D_{2;q_0}$ kills our polynomial.
We have

$$
\nabla_2 e_{1}^{m}(x_1,x_2,\dots,x_k)(x_1-x_2) =
\left(\!\!\!\begin{array}{c}
 m \\
  2 \\
\end{array}\!\!\!
\right)2k\, e_{1}^{m-2}(x_1,x_2,\dots,x_k)(x_1-x_2),
$$
while
\begin{eqnarray*}
\widetilde{D}_2 e_{1}^{m}(x_1,x_2,\dots,x_k)(x_1-x_2) & = &
\left(\left(\!\!\!\begin{array}{c}
 m \\
  3 \\
\end{array}\!\!\!
\right)6+3m(m-1)\right)e_{1}^{m-2}(x_1,x_2,\dots,x_k)(x_1-x_2)\\
 & = & (m+1) m(m-1)e_{1}^{m-2}(x_1,x_2,\dots,x_k)(x_1-x_2).
\end{eqnarray*}
Therefore
$$
D_{2;q_0}e_{1}^{m}(x_1,x_2,\dots,x_k)(x_1-x_2)=\qquad\qquad\qquad
\qquad\qquad\qquad\qquad\qquad\qquad
$$

$$
\qquad\qquad\qquad =(m(m-1)k
+q_0(m+1)m(m-1))e_{1}^{m-2}(x_1,x_2,\dots,x_k)(x_1-x_2)=0.
$$
\end{proof}
Notice that the degree of the polynomial
$e_{1}^{m}(x_1,x_2,\dots,x_k)(x_1-x_2)$ is $m+1$, hence whenever
$m+1> \left(\!\!\!\begin{array}{c}
 n \\
  2 \\
\end{array}\!\!\!
\right)$, by the previous theorem the value $q_0=-\frac{k}{m+1}$
with $2\leq k\leq n$ is singular. This shows that for each $n$,
all but finitely many of the numbers of the form $-a/b$ with $a\in
\{1,2,\dots,n\}$ and $b\in \mathbb{N}$ (the ones that show up in
Conjecture \ref{conj}) are in fact singular.

We are now in a position to proof Theorem \ref{thmmain}.

\begin{proof}[proof of Theorem \ref{thmmain}]
For every integer $d\geq 1$ and every partition $\mu$ of $d$ we
denote by $V_{\mu}$ the irreducible $\frak{S}_d$-representation
corresponding to $\mu$.

Given $m\geq 1$ and $n\geq k\geq 2$, for $1\leq i,j\leq n$, $i\neq
j$, we set
$$
p_{i,j}:=\sum_{h=1}^{n}\left(\sum_{\begin{smallmatrix}
\{i,h\}\subseteq S\subseteq
\{1,2,\dots,n\}\\
|S|=k\end{smallmatrix}}e_{1}^{m}(\mathbf{x}_{S})(x_i-x_h)-\sum_{\begin{smallmatrix}
\{j,h\}\subseteq S\subseteq
\{1,2,\dots,n\}\\
|S|=k\end{smallmatrix}}e_{1}^{m}(\mathbf{x}_{S})(x_j-x_h)\right),
$$
where $\mathbf{x}_S$ indicates the set of variables indexed by the
elements of $S$. It's easy to see that the map $p_{i,j}\mapsto
x_i-x_j$ is an isomorphism of representations of $\frak{S}_n$.
Since clearly the $p_{i,j}$'s are in the $\frak{S}_n$-module
generated by $e_{1}^{m}(x_1,x_2,\dots,x_k)(x_1-x_2)$, we have just
showed that this module contains a submodule isomorphic to
$V_{(n-1,1)}$.

All this implies the singularity of $q_0=-a/b$ with $a\in
\{1,2,\dots,n\}$, $b\in \mathbb{N}$ and $b\geq n$, since in the
Frobenius characteristic of the classical harmonics $s_{(n-1,1)}$
shows up only up to degree $n-1$. This proves the theorem.
\end{proof}
During our computer investigations we realized that we couldn't
find an example of singular value of $q_0$ which is not in the
form of Theorem \ref{thmmain}.

We risk the following conjecture.
\begin{Conj}
The numbers of the form $-a/b$ where $a,b\in \mathbb{N}$ and
$b\geq n\geq a\geq 1$ are the only singular values of $q_0$.
\end{Conj}

\section*{Acknowledgments}

We are grateful to Adriano Garsia and Nolan Wallach for useful
conversations and encouragements. The first author is pleased to
thank MPIM for partial support.

\section*{Appendix: an $f_1$ and an $f_2$ for $\partial_1 \Delta$}

\textit{Notice that in this section we use notations and results
from all the previous sections. }\newline

We want to find now a value of $c$ for which we can find an $f_2$
for $f_{1;c}$, if any exists. First we proceed as we did with our
first solution. We write
$$
f_{1;c}=\Delta^{(1)}g
$$
and
$$
g_{s,c}=a_se_{n-s-2}+b_se_{n-s-3}e_1+c_se_{n-s-4}e_{1}^2+d_se_{n-s-4}e_{2},
$$
where the coefficients are determined by the previous equations,
and of course they depend not only on $s$, but also on $n$ and
$c=c(n)$. Again we get
$$
x_1\partial_{1}^2f_{1;c}=\Delta^{(1)}\sum_{s=0}^{n-3}(s+1)s(a_{s+1}e_{n-s-3}+b_{s+1}e_{n-s-4}e_1+c_{s+1}e_{n-s-5}e_{1}^2+d_{s+1}e_{n-s-5}e_{2}),
$$
and
\begin{eqnarray*}
\sum_{j=2}^n x_j\partial_{j}^2f_{1;c} & = & 2\sum_{j=2}^n
(\partial_{j}\Delta^{(1)})x_j\partial_{j}g+\Delta^{(1)}\sum_{j=2}^n
x_j\partial_{j}^2g\\
 & = & \Delta^{(1)}\left(2\sum_{s=0}^{n-2}(P_1g_{s;c})x_{1}^s+
 \sum_{s=0}^{n-2}(\widetilde{D}_1g_{s;c})x_{1}^{s}\right).
\end{eqnarray*}

Now

\begin{eqnarray*}
 \widetilde{D}_1g_{s;c} & = & 2(n-s-3)b_se_{n-s-3}+(4(n-s-4)+2)c_se_{n-s-4}e_1\\
 & + & 2(n-s-5)d_se_{n-s-4}e_1+2(n-s-3)d_se_{n-s-3},
\end{eqnarray*}
and
\begin{eqnarray*}
 2P_1g_{s;c} & = &
 (s+2)(s+1)a_{s}e_{n-s-3}+(s+3)(s+2)b_se_{n-s-4}e_1+(n-1)(n-2)b_se_{n-s-3}\\
 & + & (s+4)(s+3)c_se_{n-s-5}e_{1}^2+2(n-1)(n-2)c_se_{n-s-4}e_1\\
 & + & (s+4)(s+3)d_se_{n-s-5}e_{2}+(n-2)(n-3)d_se_{n-s-4}e_1.
\end{eqnarray*}
Hence we can write
$$
-\widetilde{D}_1f_{1;c}=
\Delta^{(1)}\sum_{s=0}^{n-2}(\tilde{a}_s\, e_{n-s-3}+\tilde{b}_s\,
e_{n-s-4}e_1+\tilde{c}_s\, e_{n-s-5}e_{1}^2+ \tilde{d}_s\,
e_{n-s-5}e_{2})x_{1}^s,
$$
where
\begin{eqnarray*}
\tilde{a}_s & := & -(s+1)s\, a_{s+1}-
2(n-s-3)b_s-2(n-s-3)d_s-(s+2)(s+1)a_s-(n-1)(n-2)b_s\\
 & = & -(s+1)s(-1)^{s+1}\left(\frac{((-1)^{n-1}-c\,
n)}{n^2-7}
\frac{(s+1)s}{2}\frac{n-1}{n}+c\frac{(s+1)(n-1)}{2}+(s+1)(-1)^n\right)(s+2)\\
 & - & (2(n-s-3)+(n-1)(n-2))(-1)^s\left(\frac{s((-1)^{n-1}-c\,
n)}{n(n^2-7)}+ \frac{c}{2}\right)(s+2)(s+1)
\end{eqnarray*}
\begin{eqnarray*}
 & - & 2(n-s-3)(-1)^{s+1}\frac{((-1)^{n-1}-c\,
n)}{n^2-7} (s+3)(s+2)(s+1)\\
 & - & (s+2)(s+1)(-1)^s\left(\frac{((-1)^{n-1}-c\,
n)}{n^2-7}
\frac{s(s-1)}{2}\frac{n-1}{n}+c\frac{s(n-1)}{2}+s(-1)^n\right)(s+1)\\
 & = & (-1)^{s}\frac{(s+2)(s+1)}{2n(n^2-7)}\left(3(n-1)(c n +(-1)^n)s^2\begin{array}{c}
    \\
    \\
 \end{array}\right.\\
 & - & (-1)^n (21n((-1)^n -n)c+2n^2-21n+7))s\\
 & - & \left. \begin{array}{c}
    \\
    \\
 \end{array}\!\!\!\! n(n-1)(n^3+n-28)c+(-1)^n 12n(n-3)\right),
\end{eqnarray*}
\begin{eqnarray*}
\tilde{b}_s & := & -(s+1)s\, b_{s+1}-(4(n-s-4)+2)c_s -
2(n-s-5)d_s-(s+3)(s+2)b_s\\
 & - & 2(n-1)(n-2)c_s-(n-2)(n-3)d_s\\
 & = & -(s+1)s(-1)^{s+1}\left(\frac{(s+1)((-1)^{n-1}-c\,
n)}{n(n^2-7)}+ \frac{c}{2}\right)(s+3)(s+2)\\
 & - & ((4(n-s-4)+2)+2(n-1)(n-2))(-1)^s\frac{(n-1)((-1)^{n-1}-c\,
 n)}{2n(n^2-7)}(s+3)(s+2)(s+1)\\
 & - & (2(n-s-5)+(n-2)(n-3))(-1)^{s+1}\frac{((-1)^{n-1}-c\,
n)}{n^2-7} (s+3)(s+2)(s+1)\\
 & - & (s+3)(s+2)(-1)^s\left(\frac{s((-1)^{n-1}-c\,
n)}{n(n^2-7)}+ \frac{c}{2}\right)(s+2)(s+1)\\
 & = & (-1)^s \frac{(s+3)(s+2)(s+1)}{2n(n^2-7)}\left( 6(c n+(-1)^n)s+(24c n+2(-1)^n n^2+10(-1)^n)
 \right),
\end{eqnarray*}
\begin{eqnarray*}
\tilde{c}_s & := & -(s+1)s\, c_{s+1}-(s+4)(s+3)c_s\\
 & = & -(s+1)s(-1)^{s+1}\frac{(n-1)((-1)^{n-1}-c\,
 n)}{2n(n^2-7)}(s+4)(s+3)(s+2)\\
 & - & (s+4)(s+3)(-1)^s\frac{(n-1)((-1)^{n-1}-c\,
 n)}{2n(n^2-7)}(s+3)(s+2)(s+1)\\
 & = & (-1)^s\frac{(s+4)(s+3)(s+2)(s+1)}{2n(n^2-7)}(3(n-1)(c n
 +(-1)^n)),
\end{eqnarray*}
and
\begin{eqnarray*}
\tilde{d}_s & := & -(s+1)s\, d_{s+1}-(s+4)(s+3)d_s\\
 & = & -(s+1)s(-1)^{s+2}\frac{((-1)^{n-1}-c\,
n)}{n^2-7} (s+4)(s+3)(s+2)\\
 & - & (s+4)(s+3)(-1)^{s+1}\frac{((-1)^{n-1}-c\,
n)}{n^2-7} (s+3)(s+2)(s+1)\\
 & = & (-1)^s\frac{(s+4)(s+3)(s+2)(s+1)}{(n^2-7)}(-3(c n +(-1)^n))
\end{eqnarray*}

To compute $\widetilde{D}_2f_{1;c}$, first we have
$$
x_1\partial_{1}^3 f_{1;c}
=\Delta^{(1)}\sum_{s=0}^{n-4}(s+2)(s+1)s(a_{s+2}\,e_{n-4-s}+b_{s+2}\,e_{n-5-s}e_1+c_{s+2}\,e_{n-6-s}e_{1}^2+d_{s+2}\,e_{n-6-s}e_2)x_{1}^{s}.
$$

Then
\begin{eqnarray*}
\sum_{j=2}^n x_j\partial_{j}^3f_{1;c}  & = &  3\sum_{j=2}^n
(\partial_{j}^2\Delta^{(1)})x_j\partial_{j}g+ 3\sum_{j=2}^n
(\partial_{j}\Delta^{(1)})x_j\partial_{j}^2g+\Delta^{(1)}
\sum_{j=2}^n x_j\partial_{j}^3 g\\
 & = & \Delta^{(1)}\sum_{j=2}^n(6\, Q_2g_{s;c} +3\widetilde{P}_2g_{s;c}+ \widetilde{D}_2g_{s:c})x_{1}^s,
\end{eqnarray*}
where
$$
\widetilde{P}_2:=\sum_{2\leq i<j\leq
n}\frac{1}{x_i-x_j}(x_i\partial_{i}^2-x_j\partial_{j}^2).
$$
Now
\begin{eqnarray*}
6\, Q_2g_{s;c} & = & -(s+3)(s+2)(s+1)a_s
e_{n-s-4}-(s+4)(s+3)(s+2)b_se_{n-s-5}e_1\\
 & - & (s+5)(s+4)(s+3)c_s e_{n-s-6}e_{1}^2-(s+5)(s+4)(s+3)d_s
 e_{n-s-6}e_{2}\\
 & - & (n-1)(n-2)(n-3)d_s e_{n-s-4},
\end{eqnarray*}
\begin{eqnarray*}
3\, \widetilde{P}_2 g_{s;c} & = & 3(s+3)(s+2)b_se_{n-s-4}\\
 & + & 6(s+4)(s+3)c_se_{n-s-5}e_1+3(n-1)(n-2)c_se_{n-s-4}\\
 & + & 3(s+4)(s+3)d_se_{n-s-5}e_1-3(n-s-4)(n+s+1)d_se_{n-s-4},
\end{eqnarray*}
while clearly $\widetilde{D}_2g_{s;c}=0$.

Hence we have
$$
-\widetilde{D}_2f_{1;c}=\Delta^{(1)}\sum_{s=0}^{n-4}(\hat{a}_se_{n-s-4}+\hat{b}_se_{n-s-5}e_1+
\hat{c}_se_{n-s-6}e_{1}^2+\hat{d}_se_{n-s-6}e_{2})x_{1}^s,
$$
where
\begin{eqnarray*}
\hat{a}_s & := & -(s+2)(s+1)s\,
a_{s+2}+(s+3)(s+2)(s+1)a_s+(n-1)(n-2)(n-3)d_s\\
 & - & 3(s+3)(s+2)b_s-3(n-1)(n-2)c_s+3(n-s-4)(n+s+1)d_s\\
 & = & -(s+2)(s+1)s(-1)^{s+2}(s+2)\left(\frac{((-1)^{n-1}-c\,
n)}{n^2-7}
\frac{(s+1)}{2}\frac{n-1}{n}+c\frac{(n-1)}{2}+(-1)^n\right)(s+3)\\
 & + & (s+3)(s+2)(s+1)(-1)^s\left(\frac{((-1)^{n-1}-c\,
n)}{n^2-7}
\frac{s(s-1)}{2}\frac{n-1}{n}+c\frac{s(n-1)}{2}+s(-1)^n\right)(s+1)
 \end{eqnarray*}
\begin{eqnarray*}
 & + & (n-1)(n-2)(n-3)(-1)^{s+1}\frac{((-1)^{n-1}-c\,
n)}{n^2-7} (s+3)(s+2)(s+1)\\
 & - & 3(s+3)(s+2)(-1)^s\left(\frac{s((-1)^{n-1}-c\,
n)}{n(n^2-7)}+ \frac{c}{2}\right)(s+2)(s+1)\\
 & - & 3(n-1)(n-2)(-1)^s\frac{(n-1)((-1)^{n-1}-c\,
 n)}{2n(n^2-7)}(s+3)(s+2)(s+1)\\
 & + & 3(n-s-4)(n+s+1)(-1)^{s+1}\frac{((-1)^{n-1}-c\,
n)}{n^2-7} (s+3)(s+2)(s+1)\\
 & = & (-1)^s \frac{(s+3)(s+2)(s+1)}{2n(n^2-7)}\left((n-1)(c n +(-1)^n)s^2+ \begin{array}{c}
    \\
    \\
 \end{array} \right.\\
 & + & ((-n(n-1)(n^2+3n+23))c+(-1)^n(-13 n-2 n^3+9))s\\
 & + & (n(n-1)(2n^3-n^2-15n-36))c\\
 & + & \left.\begin{array}{c}
    \\
    \\
 \end{array}\!\!\!\!\!(-1)^n( -3 n^3-8 n^2-6+2 n^4-21  n)\right),
\end{eqnarray*}
\begin{eqnarray*}
\hat{b}_s & := & -(s+2)(s+1)s\, b_{s+2}+
(s+4)(s+3)(s+2)b_s-6(s+4)(s+3)c_s-3(s+4)(s+3)d_s\\
 & = & (-1)^s \frac{(s+4)(s+3)(s+2)(s+1)}{2n(n^2-7)}\left(-(6(c n +(-1)^{n}))s+18 (-1)^{1+n}+2 c n^3-32 c
 n\right),
\end{eqnarray*}
\begin{eqnarray*}
\hat{c}_s & := & -(s+2)(s+1)s\, c_{s+2}+(s+5)(s+4)(s+3)c_s\\
 & = & (-1)^s \frac{(s+5)(s+4)(s+3)(s+2)(s+1)}{2n(n^2-7)}\left(-(3(n-1))(c\, n +(-1)^n)\right)
\end{eqnarray*}
\begin{eqnarray*}
\hat{d}_s & := & -(s+2)(s+1)s\, d_{s+2}+(s+5)(s+4)(s+3)d_s\\
 & = & (-1)^s \frac{(s+5)(s+4)(s+3)(s+2)(s+1)}{n^2-7}3(c\, n +(-1)^n)
\end{eqnarray*}

Again, we set
$$
A_s:=\tilde{a}_s\, e_{n-s-3}+\tilde{b}_s\,
e_{n-s-4}e_1+\tilde{c}_s\, e_{n-s-5}e_{1}^2+ \tilde{d}_s\,
e_{n-s-5}e_{2}
$$
and
$$
B_s:=\hat{a}_se_{n-s-4}+\hat{b}_se_{n-s-5}e_1+
\hat{c}_se_{n-s-6}e_{1}^2+\hat{d}_se_{n-s-6}e_{2}.
$$

Again, we get
\begin{eqnarray} \label{semifinal2c}
\notag\frac{(-1)^{s}}{s!}\nabla_{1}^s(\nabla_2+2P_2+\nabla_{1}^2)g_{0;c}
&
= & (B_s+ \nabla_1A_s-(s+1)A_{s+1})+\qquad\qquad\\
 & - & \left(\sum_{j=0}^{s-1}(-1)^{s-1-j}\frac{j!}{s!}\nabla_{1}^{s-1-j}
(\nabla_2+2P_2+\nabla_{1}^2)A_{j}\right).
\end{eqnarray}
Now,
$$
B_s+
\nabla_1A_s-(s+1)A_{s+1}=\qquad\qquad\qquad\qquad\qquad\qquad\qquad\qquad\qquad\qquad
$$
\begin{eqnarray*}
 & = &
\hat{a}_se_{n-s-4}+\hat{b}_se_{n-s-5}e_1+
\hat{c}_se_{n-s-6}e_{1}^2+\hat{d}_se_{n-s-6}e_{2}\\
 & + & (s+3)\tilde{a}_s\, e_{n-s-4}+(s+4)\tilde{b}_s\,
e_{n-s-5}e_1+(n-1)\tilde{b}_s\, e_{n-s-4}\\
 & + & (s+5)\tilde{c}_s\,
e_{n-s-6}e_{1}^2+2(n-1)\tilde{c}_s\, e_{n-s-5}e_{1}+
(s+5)\tilde{d}_s\, e_{n-s-6}e_{2}+(n-2)\tilde{d}_s\,
e_{n-s-5}e_{1}\\
 & - & (s+1)(\tilde{a}_{s+1}\, e_{n-s-4}+\tilde{b}_{s+1}\,
e_{n-s-5}e_1+\tilde{c}_{s+1}\, e_{n-s-6}e_{1}^2+ \tilde{d}_{s+1}\,
e_{n-s-6}e_{2})\\
 & = &
 (\hat{a}_s+(s+3)\tilde{a}_s+(n-1)\tilde{b}_s-(s+1)\tilde{a}_{s+1})e_{n-s-4}\\
 & + &
 (\hat{b}_s+(s+4)\tilde{b}_s+2(n-1)\tilde{c}_s+(n-2)\tilde{d}_s-(s+1)\tilde{b}_{s+1})e_{n-s-5}e_1\\
 & + &
 (\hat{c}_s+(s+5)\tilde{c}_s-(s+1)\tilde{c}_{s+1})e_{n-s-6}e_{1}^2\\
 & + &
 (\hat{d}_s+(s+5)\tilde{d}_s-(s+1)\tilde{d}_{s+1})e_{n-s-6}e_{2}\\
 & = & (-1)^s
 \frac{(s+3)(s+2)(s+1)}{2n(n^2-7)}(-n(n-1)(n^2+3n-31)c+(-1)^n(-4n^2-17+41n-2n^3))e_{n-s-4}\\
 & + & (-1)^s
 \frac{(s+4)(s+3)(s+2)(s+1)}{2n(n^2-7)}(6(c n+(-1)^{n})s\\
  & + & (28n+2n^3)c+(-1)^n(14+4n^2))e_{n-s-5}e_1\\
  & + & (-1)^s
 \frac{(s+5)(s+4)(s+3)(s+2)(s+1)}{2n(n^2-7)}(3(n-1)(c
 n+(-1)^{n}))e_{n-s-6}e_{1}^2\\
 & + & (-1)^s
 \frac{(s+5)(s+4)(s+3)(s+2)(s+1)}{(n^2-7)}(-3(c
 n+(-1)^{n}))e_{n-s-6}e_{2}.
\end{eqnarray*}

Since
$$
(\nabla_2+2P_2+\nabla_{1}^2)(e_ke_2)=2(n-k)(n-1)e_{k-1}e_1-2ke_k,
$$
we have
$$
(\nabla_2+2P_2+\nabla_{1}^2)A_j\qquad\qquad\qquad\qquad\qquad\qquad\qquad\qquad\qquad\qquad\qquad\qquad\qquad\qquad
$$
\begin{eqnarray*}
 \qquad\quad\quad & = &
2n(j+4)\tilde{b}_je_{n-j-5}+\\
 & + & 4n(j+5)\tilde{c}_je_{n-j-6}e_1+ 2n(n-1)\tilde{c}_je_{n-j-5}\\
 & + &
 2(j+5)(n-1)\tilde{d}_je_{n-j-6}e_1-2(n-j-5)\tilde{d}_je_{n-j-5}\\
 & = & (-1)^j
 \frac{(j+4)(j+3)(j+2)(j+1)}{(n^2-7)}((3n^3-3n)c+(-1)^n(5n^2-17))e_{n-j-5}.
\end{eqnarray*}
The second term of the RHS of (\ref{semifinal2c}) is
$$
-\sum_{j=0}^{s-1}(-1)^{s-1-j}\frac{j!}{s!}\nabla_{1}^{s-1-j}
(\nabla_2+2P_2+\nabla_{1}^2)A_{j}
=\qquad\qquad\qquad\qquad\qquad\qquad
$$
\begin{eqnarray*}
& =  &
-\sum_{j=0}^{s-1}(-1)^{s-1-j}\frac{j!}{s!}\nabla_{1}^{s-1-j}
\left((-1)^j
 \frac{(j+4)(j+3)(j+2)(j+1)}{(n^2-7)}((3n^3-3n)c\right.\\
 & + & \left.\begin{array}{c}
    \\
    \\
 \end{array}\!\!\!\!(-1)^n(5n^2-17))e_{n-j-5}\right)\\
 & =  &
-\sum_{j=0}^{s-1}(-1)^{s-1-j}\frac{j!}{s!} \left((-1)^j
 \frac{(j+4)(j+3)(j+2)(j+1)}{(n^2-7)}((3n^3-3n)c\right.\\
 & + & \left.\begin{array}{c}
    \\
    \\
 \end{array}\!\!\!\!(-1)^n(5n^2-17))\frac{(s+3)!}{(j+4)!}e_{n-s-4}\right)\\
 & =  &
(-1)^s(s+3)(s+2)(s+1) \left(
 \frac{1}{(n^2-7)}((3n^3-3n)c\right.\\
 & + & \left.\begin{array}{c}
    \\
    \\
 \end{array}\!\!\!\!(-1)^n(5n^2-17))\right)\left(\sum_{j=0}^{s-1}1\right)e_{n-s-4}\\
 & =  &
(-1)^s\frac{(s+3)(s+2)(s+1)}{(n^2-7)} \left(
 (3n^3-3n)c+(-1)^n(5n^2-17)\right)s\, e_{n-s-4}.
\end{eqnarray*}

Finally, we can write (\ref{semifinal2c}) as
$$
 \nabla_{1}^s(\nabla_2+2P_2+\nabla_{1}^2)g_{0;c}=\qquad\qquad\qquad\qquad\qquad\qquad\qquad\qquad\qquad\qquad\qquad\qquad
$$
\begin{eqnarray*}
\qquad\qquad&=& \frac{(s+3)!}{2n(n^2-7)}
\left(3(n-1)(cn+(-1)^n)s^2\begin{array}{c}
   \\
   \\
\end{array}\right.\\
  & + & (n(n-1)(5n^2+3n+31)c+(-1)^n(7n+8n^3-17-4n^2))s\\
 & + & \left.\begin{array}{c}
    \\
    \\
 \end{array}\!\!\!\!\!\!\!  (-n(n-1)(n^2+17n-68))c+(-1)^n(85n-n^3-26-36n^2+2n^4)
 \right)e_{n-s-4}\\
 & + &
 \frac{(s+4)!}{2n(n^2-7)}\left(6(c n+(-1)^{n})s\begin{array}{c}
    \\
    \\
 \end{array}\!\!\!\!\!\!\! \,\, +(28n+2n^3)c+(-1)^n(14+4n^2)\right)e_{n-s-5}e_1\\
  & + &  \frac{(s+5)!}{2n(n^2-7)}(3(n-1)(c
 n+(-1)^{n}))e_{n-s-6}e_{1}^2\\
 & + &  \frac{(s+5)!}{(n^2-7)}(-3(c
 n+(-1)^{n}))e_{n-s-6}e_{2}.
\end{eqnarray*}

We reduced ourselves to solve this system of equations. We assume
that we can find a solution of the form:
\begin{equation} \label{finalf2c}
(\nabla_2+2P_2+\nabla_{1}^2)g_{0;c}=3!\alpha e_{n-4}+4! \beta
e_{n-5}e_1+5! \gamma e_{n-6}e_{1}^2+5! \delta e_{n-6}e_{2},
\end{equation}
where $\alpha,\beta,\gamma$ and $\delta$ are coefficients
depending only on $n$, and the normalization with the factorials
is made for convenience in the following computations.

We have
$$
\nabla_{1}^s\left(3!\alpha e_{n-4}+4! \beta e_{n-5}e_1+5! \gamma
e_{n-6}e_{1}^2+5! \delta e_{n-6}e_{2}\right)=
$$
\begin{eqnarray*}
 & = & \left(\alpha+\beta \, s(n-1)+\gamma\, s(s-1)(n-1)^2+\delta\,
 \frac{s(s-1)}{2}(n-1)(n-2)\right)(s+3)!e_{n-s-4}\\
 & + & (\beta+\gamma\, 2s(n-1)+\delta\,
 s(n-2))(s+4)!e_{n-s-5}e_1\\
 & + & \gamma (s+5)!e_{n-s-6}e_{1}^2+\delta (s+5)! e_{n-s-6}e_{2}.
\end{eqnarray*}

Now we have to equate the unknown coefficients to the one we have
in the system.

First we get
$$
\gamma=\frac{1}{2n(n^2-7)}(3(n-1)(c
 n+(-1)^{n}))\quad \text{and}\quad \delta= \frac{1}{(n^2-7)}(-3(c
 n+(-1)^{n})).
$$
Replacing them in the second coefficient we have
$$
\beta+\left(\frac{1}{2n(n^2-7)}(3(n-1)(c
 n+(-1)^{n}))\right)\, 2s(n-1)+\left(\frac{1}{(n^2-7)}(-3(c
 n+(-1)^{n}))\right)\,
 s(n-2)=
$$
$$
=\frac{1}{2n(n^2-7)}\left(6(c n+(-1)^{n})s\begin{array}{c}
    \\
    \\
 \end{array}\!\!\!\!\!\!\! \,\,
 +(28n+2n^3)c+(-1)^n(14+4n^2)\right),
$$
from which we get
$$
\beta=\frac{1}{n(n^2-7)}((n^3+14n)c+(-1)^n(2n^2+7)).
$$
From the other equation we get
\begin{eqnarray*}
\alpha & = &
-\left(\frac{1}{n(n^2-7)}((n^3+14n)c+(-1)^n(2n^2+7))\right)s(n-1)\\
& - & \frac{s(s-1)(n-1)^2}{2n(n^2-7)}(3(n-1)(c
 n+(-1)^{n}))\\
 & - &  \frac{s(s-1)(n-1)(n-2)}{2(n^2-7)}(-3(c
 n+(-1)^{n}))\\
 & + & \frac{1}{2n(n^2-7)}
\left(3(n-1)(cn+(-1)^n)s^2\begin{array}{c}
   \\
   \\
\end{array}\right.\\
  & + & (n(n-1)(5n^2+3n+31)c+(-1)^n(7n+8n^3-17-4n^2))s\\
 & + & \left.\begin{array}{c}
    \\
    \\
 \end{array}\!\!\!\!\!\!\!  (-n(n-1)(n^2+17n-68))c+(-1)^n(85n-n^3-26-36n^2+2n^4)
 \right)\\
 & = & \frac{1}{2n(n^2-7)}\left((3n(n-1)(n^2+n+2)c+(-1)^n(-4n+4n^3-6))s\begin{array}{c}
    \\
    \\
 \end{array}\!\!\!\!\!\!\!\right.\\
 & + & \left.\begin{array}{c}
    \\
    \\
 \end{array}\!\!\!\!\!\!\!(-n(n-1)(n^2+17n-68))c+(-1)^n(-36n^2+85n+2n^4-26-n^3)
 \right).
\end{eqnarray*}
Since we want $\alpha$ depending only on $n$, we must have
$$
3n(n-1)(n^2+n+2)c+(-1)^n(-4n+4n^3-6)=0,
$$
so
$$
c=(-1)^{n-1}\frac{2(2n^3-2n-3)}{3n(n-1)(n^2+n+2)}.
$$
We determined a value of $c$ for which we reduced all the system
to the single equation (\ref{finalf2c}).

Before computing the solution of the equation, we compute the
explicit formula for the $f_{1;c}$ for this value of $c$:
\begin{eqnarray*}
g_{s;c} & = & (-1)^{s+n}\frac{(s+1)}{6n(n^2+n+2)}(n\, s^2
+(2n^3+6n^2+15n+6)s )e_{n-s-2}\\
 & + & (-1)^{s+n}\frac{(s+2)(s+1)}{3n(n-1)(n^2+n+2)}(n\,
 s+(-2n^3+2n+3))e_{n-s-3}e_1\\
 & + &
 (-1)^{s+n}\frac{(s+3)(s+2)(s+1)}{6(n^2+n+2)}e_{n-s-4}e_{1}^2\\
 & - &
 (-1)^{s+n}\frac{(s+3)(s+2)(s+1)n}{3(n-1)(n^2+n+2)}e_{n-s-4}e_2,
\end{eqnarray*}
and from this we can write a formula for $f_1$.

Now we substitute the value of $c$ that we have found into the
coefficients:
\begin{eqnarray*}
\alpha & = & (-1)^n\frac{(6n^4+7n^3+11n^2-86n-36)}{6n(n^2+n+2)},\\
\beta & = & (-1)^n\frac{(n+2)(2n^2-4n-3)}{3n(n-1)(n^2+n+2)},\\
\gamma & = & -(-1)^n\frac{1}{2(n^2+n+2)},\\
\delta & = & (-1)^n\frac{n}{(n-1)(n^2+n+2)}.
\end{eqnarray*}
We now assume that $g_0$ is of the form
$$
u\, e_{n-3}e_1+ v\, e_{n-4}e_{1}^2+ w\, e_{n-4}e_{2}+y\,
e_{n-5}e_{1}^3+ z\, e_{n-5}e_2e_1,
$$
where $u,v,w$ and $z$ are coefficients depending on $n$ which are
to be determined.

For convenience we record the following identities:
\begin{eqnarray*}
(\nabla_2+2P_2+\nabla_{1}^2)e_{n-5}e_2e_1 & = &
10(n-1)e_{n-6}e_{1}^2+10ne_{n-6}e_2+(2n(n-2)-2(n-5))e_{n-5}e_1;\\
(\nabla_2+2P_2+\nabla_{1}^2)e_{n-5}e_{1}^3 & = & 30n
e_{n-6}e_{1}^2+6n(n-1)e_{n-5}e_1;\\
(\nabla_2+2P_2+\nabla_{1}^2)e_{n-4}e_{2} & = &
8(n-1)e_{n-5}e_1-2(n-4)e_{n-4};\\
(\nabla_2+2P_2+\nabla_{1}^2)e_{n-4}e_{1}^2 & = &
16ne_{n-5}e_1+2n(n-1)e_{n-4};\\
(\nabla_2+2P_2+\nabla_{1}^2)e_{n-3}e_{1} & = & 6ne_{n-4}.
\end{eqnarray*}
we get
\begin{eqnarray*}
(\nabla_2+2P_2+\nabla_{1}^2)g_0 & = & 10n\, z\,
e_{n-6}e_2+\left(10(n-1)\, z+30n\, y \right)e_{n-6}e_{1}^2\\
 & + & \left((2n(n-2)-2(n-5))z+6n(n-1)y+8(n-1)w+16n\,
 v\right)e_{n-5}e_1\\
  & + & \left(-2(n-4)w+2n(n-1)v+6n\, u\right)e_{n-4}.
\end{eqnarray*}
Equating coefficients we have:
\begin{eqnarray*}
z & = & (-1)^n\frac{-6}{n(n^2+n+2)};\\
y & = & (-1)^n\frac{-2}{n^2(n^2+n+2)};\\
w & = &
-\frac{1}{2(n^2+n+2)(n^2-7)}\left((-12n^2(n^2+n+2))u+(12n^4+7n^3+31n^2-168n-48)\right);\\
v & = &\frac{1}{2}
\frac{1}{(n-1)(n^2-7)(n^2+n+2)n^2}\left((-6n^2(n^2+n+2)(n-1)^2)u\right.\\
 & + & \left.(6n^6-5n^5+10n^4-138n^3+179n^2-22n+60)\right),
\end{eqnarray*}
and $u$ is arbitrary. Hence we got a family $g_{0;u}$ of
solutions.

\end{document}